# The large-scale charging scheduling problem for fleet batteries: Lagrangian decomposition with time-block reformulations

Sunney Fotedar[a], Jiaming Wu[a], Balázs Kulcsár[b] and Rebecka Jörnsten[c]

[a]*Department of Architecture and Civil Engineering, Chalmers University of Technology, Gothenburg, Sweden*
[b]*Department of Electrical Engineering, Chalmers University of Technology, Gothenburg, Sweden*
[c]*Department of Mathematical Sciences, Chalmers University of Technology, Gothenburg, Sweden*

ABSTRACT

There is a rise in the need for efficient battery charging methods due to high penetration of electromobility solutions. Battery swapping, a technique in which fully or partially depleted batteries are exchanged and then transported to a central facility for charging, introduces a unique scheduling problem. For scenarios involving a large number of batteries (common in micromobility and energy storage applications), commercial solvers and existing methods do not yield optimal or near-optimal solutions in reasonable time due to high computational complexity. Our study presents a novel approach that combines variable layering with Lagrangian decomposition. We develop a new tighter time-block reformulation for one of the Lagrangian sub-problems, enhancing convergence rates when used with our partial-variable fixing Lagrangian heuristic. We also propose an ergodic-iterate-based local search method to further improve the solution quality. Lower bounds are improved by learning the relation between Lagrangian multipliers and electricity cost. Our extensive benchmarks show superior computational performance against a commercial mixed-integer linear programming solver (Gurobi), a constraint-programming satisfiability solver (Google's CP-SAT). We achieved, on average, a 43% lower objective value compared to state-of-the-art methods. In 71% of the instances, we obtained near-optimal solutions (optimality gap less than 6%), and 93% of the instances were below 10%. We obtained feasible solutions for all instances, compared to only 65% feasibility using incumbent methods. The developed exact method aims to support future research on charging scheduling, especially important for micromobility industry, vehicle-to-grid (V2G) applications, and second-life utilization of batteries. Furthermore, the developed polyhedral insights can be useful in other scheduling problems with a common underlying mathematical structure.

## 1. Introduction

Electromobility plays an important role in reducing carbon emissions and breaking dependency on fossil fuels. As the number of electric vehicles grows, complemented by the expansion of commercial charging infrastructure and the adoption of battery swapping (MIT), there emerges a critical need for comprehensive analysis of the underlying charging scheduling problem. This necessity becomes apparent in ensuring the efficient management of charging resources to meet the evolving demands of an electrified transportation ecosystem. In this work, we use micromobility as an example to present a novel approach in addressing this critical and challenging problem. This use case is particularly important as it involves the charging of a large volume of depleted batteries, hence necessitating advanced optimization algorithms. Moreover, it can be easily generalized into other contexts such as vehicle-to-grid charging scheduling in large charging/parking facilities, drone charging and the second-life utilization of batteries for energy storage applications (Dong et al., 2023). While most existing efforts have focused on routing problems for electric vehicles, there remains a general lack of understanding of how to optimize charging schedules in a mathematically robust way—despite the growing commercial interest in this area.

Micromobility is rapidly becoming a key solution to the pervasive first/last mile mobility problems, with e-scooters and e-bikes quickly dominating the urban transportation landscape (Zhao et al., 2024). A total of 157 million trips were made in the US and Canada in 2023 (NACTO, 2024). Most large cities are flooded with shared micromobility services such as (but not limited to) e-scooters and e-bikes. The micromobility market has been growing rapidly but there is a need for technologies and methods to improve operational efficiency and also meet tighter regulations imposed by local authorities. For example, authorities in Stockholm have implemented regulations to limit e-scooter fleet sizes and enforce parking restrictions (Bloomberg, 2022). It is now a common practice that municipalities issue tenders where operators compete based on their performance. Consequently, companies place significant emphasis on technologies that enhance operational efficiency to secure a stronger market presence. Swappable battery is one such example that has been deployed in several services of dockless shared e-scooters, which remarkably eases the process of charging depleted vehicles and thus has become the predominant approach (Gauquelin, 2020; Leurent, 2022; Zhan et al., 2022). However, studies have indicated that approximately 1.5 batteries per vehicle are typically required for this approach (Link, 2020). In some cases, based on communication with local operators in Sweden, each e-scooter may be equipped with up to three batteries at the time of purchase. With hundreds or even thousands of e-scooters operating in a city, operators typically face a challenging battery charging problem that must be dealt with on a daily basis. To our knowledge, the problem of optimizing large volume of battery charging at the warehouse is largely overlooked and, as will be shown, is far from trivial.

✉ sunneyfotedar@gmail.com (S. Fotedar); jiaming.wu@chalmers.se (J. Wu); kulcsar@chalmers.se (B. Kulcsár); jornsten@chalmers.se (R. Jörnsten)

In this paper, we address the charging scheduling problem for large fleets. It is relevant for most operators who are concerned with: a) costs of electricity, b) labor costs required to implement the developed schedule, and c) batched demands on a daily basis d) high investment cost for buying additional batteries. The price of electricity, often dynamic in many cities, underlines the necessity and complexity of a charging schedule. Labor factors must be considered since the industry currently relies on their staffs to manually plug/unplug the batteries. Furthermore, fully charged batteries are needed several batches a day in the current practice in coordination with shifts of field staff. The above-described charging problem presents a significant operational opportunity but also poses a challenging mathematical problem: solving large-scale scheduling problem to near-optimality. Apart from providing good feasible solutions, it is crucial for such operational models to provide lower bounds or guarantees on the optimality gap. Our tests indicate that current state-of-the-art approaches often fail to yield high-quality solutions and, in some cases, do not produce any feasible solutions due to the large dimension of the problems. This limitation hampers the potential benefits of employing a centralized charging schedule, which is resolved in the present research.

### 1.1. Literature Review

This section first provides an overview of the literature on the battery charging scheduling problem and subsequently discusses relevant solution methodologies.

**Battery charging scheduling**

To the best of our knowledge, the Battery Charging Station Scheduling (BCSS) problem has primarily been studied from two perspectives: vehicle or service-oriented approaches and grid or power-oriented approaches.

From the vehicle or service perspective, the focus is typically on optimizing charging or battery swapping operations for electric vehicles (EVs), often in conjunction with routing decisions (Tan et al., 2019; Cui et al., 2023; Nayak and Misra, 2024; Pelletier et al., 2018; Lai and Li, 2024). Many studies address joint routing and scheduling problems, including intra-route charging or swapping for commercial EVs (Raeesi and Zografos, 2022). For example, Adler and Mirchandani (2014) aim to minimize EV wait times at battery swapping stations through efficient routing. Coordination between EV operations and battery swapping schedules has also been explored in multi-stakeholder settings, such as integrating an isolated microgrid (IMG) with a battery swapping station using a bi-level optimal scheduling model (Li et al., 2018).

From the grid and power systems perspective, research often centers on the role of centralized charging stations (CCS) in energy systems (Sun et al., 2018; Tan et al., 2019; Nurre et al., 2014). In the absence of electricity market bidding, these operations must contend with fluctuating electricity prices. CCSs can contribute to grid stability through peak shaving and also enable operators to sell stored energy to microgrids or local flexibility markets. This line of research emphasizes the integration of charging stations into broader power systems and energy markets, rather than focusing solely on vehicle-level decisions.

Our model is closest to the scheduling model in Tan et al. (2019) and Raviv (2012), although both are still for EVs with several differences. To the best of our knowledge, our model is different in several ways as compared to existing literature and has resemblance to common practice in micro-mobility businesses (Gössling, 2020; VOI Technology, 2024). Furthermore, we believe that this is the first thorough polyhedral investigation of large-scale charging scheduling problems for central charging facility of an e-scooter fleet operator. Following are some of the key difference to existing studies:

- Fixed charging rates: As opposed to EV charging (Tan et al., 2019; Raeesi and Zografos, 2022; Pelletier et al., 2018), it is not common to have continuously-controllable charging rates for e-scooter batteries. In EV charging scheduling literature continuous variables are used to model charging variables but our proposed model is a pure binary linear programming model. This essentially differentiates the polyhedral structure of the problem significantly and makes it less tractable. Hence, a specialized approach is necessary.

- Manual loading/unloading (*switching*): Our model includes several knapsack constraints to limit the number of loading/unloading (henceforth, switching) in each time period. In contrast to EVs, e-scooter batteries are smaller, less expensive, and generally subjected to less intensive use cycles. Therefore, the financial and environmental incentives to implement smart charging for e-scooters may be perceived as less significant/practical. Hence, such constraints are necessary and computationally challenging.

- Objective function: In Raviv (2012) and Tan et al. (2019), there is no cost for switching. In Raviv (2012), it is allowed to not meet demand for charging and that is penalized in the objective function, where as in Tan et al. (2019) only time-of-use electricity price is considered. Contrastingly, we consider a *lexicographic minimization* (Ehrgott, 2005, Ch. 5.1) of time-of-use electricity cost and total number of switching (i.e. loading/unloading).

- Size of the instance: In Tan et al. (2019) the instances can involve up to 200 batteries, with no limitations on the number of switchings, which necessitates several knapsack constraints in our model as well. In our test instances, we handle up to 400 batteries. Additionally, even for instances with an equal number of batteries, our model incorporates significantly more binary variables.

- Demand window: Each depleted battery (DB) in stock is associated with a demand window, which refers to the specific time of day by which the aforementioned DB should be charged to the desired level, removed from the charging bay, and prepared for dispatch.

Our contributions exploit the polyhedral structure of the problem and can be effectively utilized in various charging scheduling problems characterized by similar mathematical properties. The proposed model and the developed exact method are general and can be applied to a variety of use cases directly or adapted with some modifications. These include, but are not limited to, Vehicle-to-Grid (V2G) applications and the second-life utilization of batteries for energy storage solutions.

**Solution methodologies for large-scale structured optimization problems**

Many practical applications of combinatorial, and large-scale integer linear optimization consist of highly structured and solvable sub-problems, linked by complicating constraints or variables. This has driven the development of classic decomposition frameworks such as Benders' decomposition (Benders, 1962), Brand and Price (Barnhart et al., 1998), and Subgradient optimization (Held and Karp, 1970). Over time, these methods have seen various enhancements. We have concentrated on Subgradient optimization, attributing to its block decomposability and our success in exploiting polyhedral structure of subproblems. Metaheuristic and other inexact approaches (Coelho and Vanhoucke, 2023; Ki et al., 2018; Özarık et al., 2023) are also popular in the literature. However, the focus of this work is to develop an exact approach with provable guarantees on the optimality gap and can be applied to future extensions of the problems as well. That said, we also incorporate a heuristic method to enhance solution quality where appropriate. One of the main reasons for choosing the subgradient optimization framework is its ease of implementation and its surprisingly strong performance on our problem instances. In fact, initial vanilla implementations of subgradient optimization outperformed Benders' decomposition and other exact methods. As a result, our discussion will primarily center on Lagrangian relaxation and subgradient optimization, in order to maintain technical focus and clarity.

*Background*: The essence of the Lagrangian dual approach, particularly Lagrangian relaxation, lies in the relaxation of a subset of constraints, thereby simplifying the original problem into more manageable independent Lagrangian sub-problems. This approach proves beneficial, especially within integer programming, where the absence of strong duality necessitates careful consideration of the formulation's *compactness* and *tightness*. A critical aspect involves the selection of constraints to relax, as improper selection can lead to weak lower bounds, equivalent to those obtained by a mere relaxation of integrality. In addition to selecting appropriate constraints for relaxation, it is vital to determine an effective method for computing the subgradient, feasible solution, and a suitable step length parameter. These components are crucial for the success of the method. Without loss of generality, given Lagrangian multiplier $\boldsymbol{\mu}$, the Lagrangian dual problem is defined as $h_D^* = \max_{\boldsymbol{\mu}} h(\boldsymbol{\mu})$, where $h(\boldsymbol{\mu}) := \min_{\mathbf{x} \in X} L(\mathbf{x}; \boldsymbol{\mu})$ is the dual function and $L(\mathbf{x}; \boldsymbol{\mu})$ is an affine Lagrangian function of $\mathbf{x}$ and finally the set $X$ is polyhedral in this study. It is well known that $h(\boldsymbol{\mu})$ is a piece-wise linear concave function and therefore, the dual problem is a non-smooth convex optimization problem (Wolsey, 2001, Thm. 10.3). Some other popular techniques for solving large-scale optimization problems are combining augmented Lagrangian method with Alternating Direction Method of Multipliers (ADMM) (Boyd et al., 2011), stabilized cutting-plane methods (Bonnans et al., 2006) and some other non-smooth convex optimization techniques are summarized in (Sagastizábal, 2012).

The function $h(\boldsymbol{\mu})$ evaluated at any point $\boldsymbol{\mu}$ may result in decomposition of the feasible set $X$ into independent Lagrangian subproblems. After computing the solutions to the Lagrangian subproblems, the subgradient is obtained at no additional cost, $\mathbf{g}^i \in \partial h(\boldsymbol{\mu})$, where $\partial h(\boldsymbol{\mu})$ represents the subdifferential at $\boldsymbol{\mu}$, and $i$ denotes the iteration counter. Hence, we have a non-zero direction to move whenever we have primal infeasibility. Although standard implementation of subgradient optimization is easy but the downside is that such implementation suffer from so-called *stabalization issues*. Techniques such as bundle method (more recent proximal bundle subgradient (Kiwiel, 2006)) and modified deflected subgradient (MDS) (Belgacem and Amir, 2018) (combining the classic modified gradient technique (Camerini et al., 1975) and average direction strategy (Sherali and Ulular, 1989)) utilize information from previous iterations to stabilize dual oscillations and these methods also preserve convergence properties albeit at the cost of additional computational effort especially in the case of bundle methods. In modified deflected subgradient (MDS), the main idea is to prevent the oscillating behavior of the dual value by weighted average of a current subgradient and the previous step direction. Let $\mathbf{g}^i$ be the current subgradient of $h$ at $\boldsymbol{\mu}^i$ and $\mathbf{d}^i$ be the step direction then we employ the formulae

$$\mathbf{d}^i := \mathbf{g}^i + Y_{\text{MDS}}^i \mathbf{d}^{i-1}, \qquad Y_{\text{MDS}}^i := (1-\alpha_i) Y_{\text{MGT}}^i + \alpha_i Y_{\text{ADS}}^i, \qquad (1a)$$

$$Y_{\text{MGT}}^i := \left[\zeta_i \frac{\langle \mathbf{g}^i, \mathbf{d}^{i-1} \rangle}{\|\mathbf{d}^{i-1}\|^2}\right]_+, \qquad Y_{\text{ADS}}^i := \frac{\|\mathbf{g}^i\|}{\|\mathbf{d}^{i-1}\|}, \qquad (1b)$$

where, as per (Belgacem and Amir, 2018), we let $\zeta_i = \frac{1}{2-\alpha_i}$ and $\alpha_i = \left[\frac{-\langle \mathbf{g}^i, \mathbf{d}^{i-1} \rangle}{\|\mathbf{g}^i\|\|\mathbf{d}^{i-1}\|}\right]_+$. In this work, the updated multiplier is defined as $\boldsymbol{\mu}^{i+1} = \boldsymbol{\mu}^i + \tau \mathbf{d}^i$, where the step length, $\tau$, primarily follows the classic rule proposed by Polyak (Polyak, 1969). Specifically, $\tau$ is calculated using the formula $\tau = \theta^i \frac{\bar{h} - h(\boldsymbol{\mu}^i)}{\|\mathbf{d}^i\|^2}$, where $\bar{h}$ represents the best known upper bound, i.e., the objective value and $h(\boldsymbol{\mu}^i)$ is the current dual value at iteration $i$. The step length parameter, $\theta^i > 0$, is reduced if no improvement in the lower bounds is identified.

### 1.2. Contributions

In this work, we address a large-scale battery charging scheduling problem, applying polyhedral theory and convex non-smooth optimization techniques. The problem studied is highly relevant not only to the micromobility industry, which serves as an example in this paper, but also to a wide range of charging scheduling problems across the broader electromobility sector. Our principal contributions are outlined as follows:

- We introduce a new optimization model for the battery charging scheduling problem that can accommodate a range of realistic operational constraints.

- We introduce a variable-layering technique to enable efficient decomposition and lagrangian relaxation of the original problem.

- We propose a novel time-block reformulation for one of the Lagrangian subproblems to tackle the computational challenges inherent in simply using the original formulation. We prove that the new time-block formulation is *tight*, meaning the feasible set it defines is a strict subset of the polyhedron defined by the original formulation, while still maintaining the same optimal value/solution.

- We present a new partial variable-fixing heuristic that yields feasible solutions and corresponding upper bounds. This heuristic is based on partially fixing variables derived from the subproblem solutions and then applying a new, compact disjunctive time-block bin-packing formulation to expedite convergence.

- Finally, we present a sensitivity analysis to investigate the impacts of various constraints and pricing conditions on the problem.

We compare our approach with alternative methods, including solvers such as Gurobi and Google's CP-SAT. Variations of our proposed algorithm are also tested to assess the sensitivity of its performance to changes in different algorithmic components.

### 1.3. Outline

In Section 2, we formally define the optimization problem. Furthermore, in Section 3 we present a binary optimization model. We provide an overview of the problem's dimension and examine how solution quality is affected as instance size increases. In Section 4, we propose introducing copies of a subset of variables of the optimization model and relaxing certain constraints, thereby obtaining a Lagrangian decomposition into two independent subproblems. We identify important polyhedral properties associated with these subproblems. In Section 4.1, we reformulate a subproblem that had been a computational bottleneck, resulting in substantial reductions in computational time and enabling more iterations of the subgradient optimization method within the allotted time frame. Section 5 outlines various modifications to standard subgradient optimization method including a Lagrangian heuristic (Section 5.1) and a local search strategy (Section 5.2). The paper concludes with Section 6, which includes computational benchmarks against solvers like CP-SAT and Gurobi, and a comprehensive analysis of the effects of different subgradient optimization components on the optimality gap, lower bounds, upper bounds, and primal-dual integrals.

## 2. Problem statement

Consider the role of an operator overseeing a large fleet of e-scooters and a charging facility (such as a warehouse). At the end of each day, depleted batteries from deployed e-scooters across the city are collected and sent back to the charging facility, and each battery's state of charge (SOC) is known. The operator also has a good estimation of next day's demand of batteries from experience and thus needs to schedule charging accordingly. Note that, fully charged batteries are often requested in batches to align with the shifts of field battery-swapping staff, which is also the setup of the present research.

The operator's task involves efficiently managing these depleted batteries (DBs) within a set planning horizon of $T$ hours. Upon arrival at the charging facility, both incoming and existing DBs are evaluated based on their SOC and assigned to one of several time-specific demand windows, denoted as $\mathcal{L} = 1, \ldots, L$. The objective is to ensure that batteries allocated to demand windows from 1 to $L-1$ are fully charged by the end of each respective period. While the times mentioned here — 14:00, 19:00, and 23:59 — serve as examples in Figure 1, operators can adjust these times based on their specific operational needs. Each window ends at a predefined time by which the charged batteries must be boxed and ready for deployment.

The final demand window, ending just before midnight (in the example), has a unique requirement: batteries assigned to this window must be charged to at least $\alpha$ SOC, after which they can be rolled over to the next day's operations. The charging of this last batch is primarily to mitigate battery degradation which is a detrimental result of leaving batteries at low SOC levels. The whole operation is illustrated in Figure 1, where batteries arrive at the facility, are charged at designated bays, and are then prepared for dispatch at specified times.

The main goal of the operation is to minimize time-of-use electricity costs, which fluctuate hourly. Additionally, it is favorable to limit and minimize manual efforts in plugging and unplugging batteries, thereby preventing excessive labor workload and optimizing operational efficiency. In Problem statement 1 (relevant notations in Table 1), we formally define the so-called *Discrete-Charging Scheduling Problem*.

**Problem Statement 1** (Discrete-Charging Scheduling Problem (D-CSP)). *Consider a set of batteries, denoted by $\mathcal{B}_\ell$, assigned to a unique demand window $\ell \in \mathcal{L}$, where $\mathcal{B} = \bigcup_{\ell \in \mathcal{L}} \mathcal{B}_\ell$ represents the total set of batteries. Additionally, there is an identical set of charging bays/ports, $\mathcal{N}$. For each battery $j \in \mathcal{B}_\ell$, a specific charging requirement $p_j \in \mathbb{Z}_+$ (measured in hours) must be met within its demand window, i.e in time periods $[1, \ldots, n_\ell]$, where $n_\ell \leq T$ and $T$ represents the total available time periods. The charging can be paused (allowing for preemption) to capitalize on varying electricity prices throughout the day. A restriction is placed on the total number of battery switchings (loading/unloading) at $\gamma \in \mathbb{Z}_+$ for each time period. The goal is to minimize the electricity costs, based on time-of-use, and the total number of battery switchings (loading/unloading), prioritized in that sequence, across all possible schedules referred to as* lexicographic minimization *of the two objectives.*

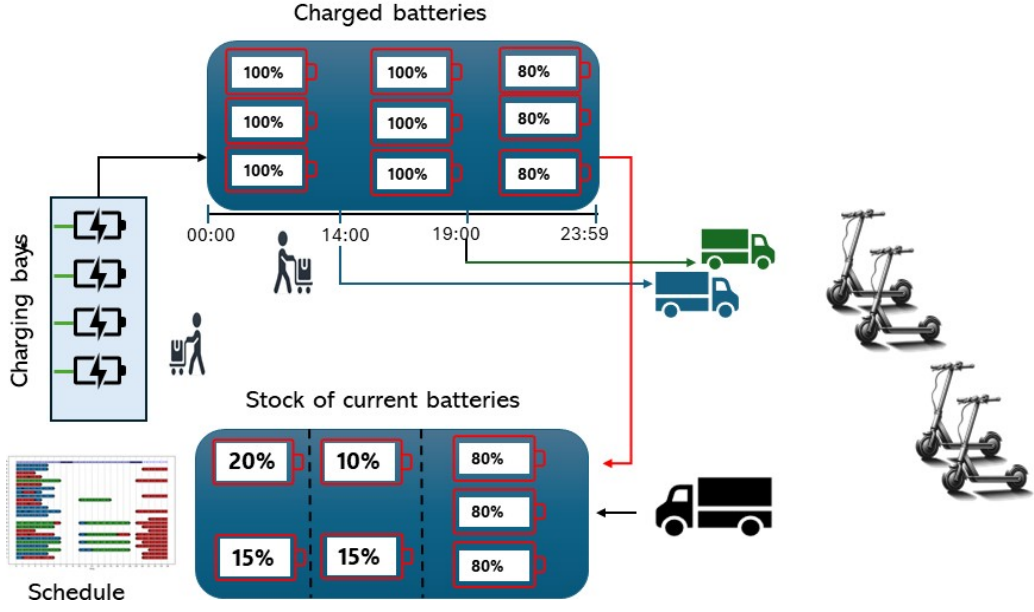

**Figure 1:** An overview of the D-CSP with three demand windows.

**Table 1**
Notation used for the D-CSP

| Sets | Description |
| --- | --- |
| $\mathcal{B} := \{1, \ldots, B\}$ | set of depleted batteries (DBs) in the warehouse |
| $\mathcal{N} := \{1, \ldots, K\}$ | set of identical charging bays (CBs) or ports |
| $\mathcal{L} := \{1, \ldots, L\}$ | set of demand windows (DW) |
| $\mathcal{T} := \{1, \ldots, T\}$ | set of time periods |
| $\mathcal{T}_\ell := \{1, \ldots, n_\ell\}$ | set of time periods associated with each DW $\ell \in \mathcal{L}$. Note that $\mathcal{T}_L = \mathcal{T}$ |
| $\mathcal{B}_\ell \subseteq \mathcal{B}$ | set of DBs to be delivered by the end of DW $\ell \in \mathcal{L}$. Note that $\underset{\ell \in \mathcal{L}}{\cup} \mathcal{B}_\ell = \mathcal{B}$ and $\underset{\ell \in \mathcal{L}}{\cap} \mathcal{B}_\ell = \emptyset$ |
| **Variables** | **Description** |
| $x_{jkt} \in \{0,1\}$ | $= 1$ if battery $j \in \mathcal{B}_\ell$ is assigned to port $k \in \mathcal{N}$ during time period $t \in \mathcal{T}_\ell, \ell \in \mathcal{L}$ |
| $y_{kt} \in \{0,1\}$ | $= 1$ if any battery $j \in \mathcal{B}$ is either loaded or offloaded from port $k \in \mathcal{N}$ after time period $t \in \mathcal{T} \setminus \{T\}$ |
| **Parameters** | **Description** |
| $n_\ell \in \mathcal{T}$ | is the last time period, by which the batteries in $\mathcal{B}_\ell, \ell \in \mathcal{L}$ must be delivered. Note that $n_\ell \leq T$ |
| $p_j \in \mathbb{Z}_+$ | total hours needed to fully charge battery $j \in \mathcal{B}_\ell, \ell \in \mathcal{L}$ |
| $c_t \in \mathbb{Q}_+$ | electricity price associated with each time period $t \in \mathcal{T}$ |
| $\alpha \in (0, 1]$ | fraction of battery to be charged if a battery belongs to the last demand window $\ell = L$ |
| $\gamma \in \mathbb{Z}_+$ | upper limit on the number of switchings at the end of each time period |
| $w \in \mathbb{R}_+$ | $:= \frac{\min_{t \in \mathcal{T}}\{c_t\}}{((T-1) \times \gamma) + 1}$ is the augmentation parameter used for lexicographic minimization |

## 3. Mathematical model

We describe a binary linear programming model for the D-CSP. The notations used are summarized in Table. 1.

$$z^* := \min \left( \sum_{\ell \in \mathcal{L}} \sum_{t \in \mathcal{T}_\ell} c_t \left( \sum_{j \in \mathcal{B}_\ell} \sum_{k \in \mathcal{N}} x_{jkt} \right) + w \sum_{k \in \mathcal{N}} \sum_{t \in \mathcal{T} \setminus T} y_{kt} \right), \tag{2a}$$

$$\text{s.t.} \sum_{t \in \mathcal{T}_L} \left( \sum_{k \in \mathcal{N}} x_{jkt} \right) \geq \lceil \alpha \, p_j \rceil, \qquad j \in \mathcal{B}_L, \qquad (2b)$$

$$\sum_{t \in \mathcal{T}_\ell} \left( \sum_{k \in \mathcal{N}} x_{jkt} \right) = p_j, \qquad j \in \mathcal{B}_\ell, \ell \in \mathcal{L} \setminus \{L\}, \qquad (2c)$$

$$\sum_{\ell \in \mathcal{L}} \sum_{j \in \mathcal{B}_\ell} x_{jkt} \leq 1, \qquad k \in \mathcal{N}, t \in \mathcal{T}, \qquad (2d)$$

$$\sum_{k \in \mathcal{N}} x_{jkt} \leq 1, \qquad j \in \mathcal{B}_\ell, t \in \mathcal{T}_\ell, \ell \in \mathcal{L}, \qquad (2e)$$

$$y_{kt} \geq \left( \sum_{\ell \in \mathcal{L}} \sum_{j \in \mathcal{B}_\ell} x_{jkt+1} - \sum_{\ell \in \mathcal{L}} \sum_{j \in \mathcal{B}_\ell} x_{jkt} \right), \qquad k \in \mathcal{N}, t \in \mathcal{T} \setminus \{T\}, \qquad (2f)$$

$$y_{kt} \geq (x_{jkt} - x_{jkt+1}), \qquad j \in \mathcal{B}_\ell, k \in \mathcal{N}, t \in \mathcal{T}_\ell \setminus \{T\}, \ell \in \mathcal{L}, \qquad (2g)$$

$$\sum_{k \in \mathcal{N}} y_{kt} \leq \gamma, \qquad t \in \mathcal{T} \setminus \{T\}, \qquad (2h)$$

$$x_{jkt} \in \{0, 1\}, \qquad j \in \mathcal{B}_\ell, k \in \mathcal{N}, t \in \mathcal{T}_\ell, \ell \in \mathcal{L}, \qquad (2i)$$

$$y_{kt} \in \{0, 1\}, \qquad k \in \mathcal{N}, t \in \mathcal{T} \setminus \{T\}. \qquad (2j)$$

The objective (2a) is to minimize lexicographically two objective functions, i.e., time-of-use electricity cost and total number of switchings. The augmentation parameter $w$ is in accordance to Mavrotas (2009) to ensure that first we minimize electricity cost and then number of switchings (cf. Table 1). The decision variable $x_{jkt}$ governs whether a battery $j$ is being charged via charger $k$ in the discrete time period $t$, while $y_{kt}$ indicates whether the usage of charger $k$ has changed at $t$, which is referred to as a switch.

The first set of constraints (2b) ensures that all the batteries assigned to the last demand window, i.e., $\mathcal{B}_L$, are at least charged $\lceil \alpha \, p_j \rceil$ hours. This requirement is to maintain good battery health by moderating battery state-of-charge (Severson et al., 2019). Furthermore, (2c) apply restriction that the demand has to be met for each demand window except the last. Constraint (2d) ensure that at most one battery can be assigned to a port $k \in \mathcal{N}$ in a time period $t \in \mathcal{T}$. Constraint (2e) ensure that the same battery $j$ is not allocated to more than one port in the same time period $t \in \mathcal{T}$. The next two constraints (2f)–(2g) ensure that we set $y_{kt} = 1$ whenever we load a battery on a port or unload it after time period $t \in \mathcal{T} \setminus \{T\}$ on port $k \in \mathcal{N}$. Constraint (2f) specifically ensure that if no battery is charged on a given port $k \in \mathcal{N}$ and then in the upcoming time period if one of the batteries is charged at port $k$ then we impose a constraint $y_{kt} \geq 1$. To address the case when both unloading of battery $j$ and loading of a distinct battery $j'$ happens at the end of time period $t$, (2g) imposes $y_{kt} \geq 1$. The last constraint (2h) restricts the total number of switchings possible at the end of any time period $t \in \mathcal{T} \setminus \{T\}$. We also note that the constraints (2f), (2g) link consecutive time periods which makes $y_{kt}$ a complicating variable and complicates the process of solving the problem as discussed in upcoming sections. Finally, all the variables are restricted to be binary.

### 3.1. Problem dimension and model structure

We offer a concise overview of the problem's dimension and how it grows with the increasing number of batteries and charging bays (CBs). This discussion aims to underscore the necessity of decomposing the problem due to the high computational intensity. We set $T = 24$ for day-ahead scheduling through hourly discretization. As the quantities of batteries and CBs grow, we observe an increase in possible configurations, as depicted by the number of variables and constraints in Table 2. The table's first column enumerates the number of batteries ($|\mathcal{B}|$), the second the number of CBs ($|\mathcal{N}|$), and the third is the limitations on total switchings per time period ($\gamma$). The next columns detail the count of variables and constraints, respectively. The final four columns display the best known estimate of the optimal value ($\bar{z}^*$), the best known lower bound on the optimal value ($\underline{z}^*$), the total solution time, and the optimality gap, respectively. For these calculations, we utilize a state-of-the-art commercial MIP solver, Gurobi.

Further computational setup and specifics will be discussed later in Section 6. These preliminary findings clearly show that, for larger instances involving 200-400 batteries, the optimality gap is significantly large indicating that much better solutions can exist. Notably, for the instances with 350 and 400 batteries, Gurobi failed to obtain any feasible solution in two hours. The relative dimensions of $\mathcal{B}$, $\mathcal{N}$, and $\gamma$ are elaborated upon in the computational section, which also reflects on the quality of solutions. An extensive usage of RAM is also observed along the expansion of the branch-and-bound tree, typically a built in method in prevailing commercial solvers.

**Table 2**
Increase in the computational effort when the dimension of the problem is increased and computational experiments with Gurobi.

| $|\mathcal{B}|$ | $|\mathcal{N}|$ | $\gamma$ | # of vars. | # of consts. | $\bar{z}^*$ | $\underline{z}^*$ | Time[sec.] | Gap[%] |
|---|---|---|---|---|---|---|---|---|
| 100 | 50 | 13 | 91150 | 119274 | 4700.32 | 4259.74 | 3600 | 9.37 |
| 200 | 100 | 25 | 362300 | 468524 | 11482.23 | 7426.75 | 3600 | 35.32 |
| 300 | 140 | 35 | 759220 | 978304 | 17716.19 | 11374.26 | 3600 | 35.79 |
| 350 | 140 | 35 | 885220 | 1140254 | - | 13718.49 | 7200* | 100 |
| 400 | 140 | 45 | 1011220 | 1302204 | - | 16296.94 | 7200* | 100 |

The model structure for very large-scale optimization problems can often be exploited by utilizing different decomposition and reformulation techniques. The first question that must be investigated for any binary or mixed integer linear programming formulation like model (2) is if there exists an alternative formulation with no integrality gap, i.e. so-called *perfect formulation*. Perfect formulations are available for certain classical combinatorial optimization problems, such as assignment, shortest path, maximum flow, min-cost network flow, bipartite matching and so on, but generally most integer programming problems do not have such perfect formulations available and requires generation of exponentially many cuts. Such problems are also referred to as $\mathcal{NP}$-hard problems. We provide a formal proof of $\mathcal{NP}$-hardness by problem reduction later as it is almost trivial to reduce a simple problem to $\mathcal{NP}$-hard problem but the reverse is not possible.

Table 2 demonstrates that attempts to derive a tighter formulation by adding classes of valid inequalities (VIs) only increases the problem size, further straining the computer's RAM, which already struggles with instances of more than 300 batteries. Therefore, decomposition techniques, particularly Lagrangian relaxation, might be more viable. This approach allows us to break the problem into smaller, independent Lagrangian sub-problems that can be solved in parallel. Decomposing the original model (2) into manageable subproblems is, however, not straightforward. We briefly evaluate three standard options before presenting our solution:

1. Decompose over time $\mathcal{T}$: it involves relaxing inter-temporal constraints that link different periods, such as (2b) and (2c), as well as (2f) to (2g). Relaxing these constraints results in a Lagrangian multiplier vector of dimension $BL + NT + BNT$. For example, with $B = 300$, $N = 140$, $T = 24$, and $L = 3$, the vector dimension reaches 1,012,260. However, as Table 2 indicates, with already 759,220 variables in the original formulation, this approach does not instill confidence in our ability to solve the problem effectively. Most importantly, the remaining sub-problem with an affine objective function separable over **x** and **y** including constraints (2d)-(2e) for the problem in **x** and (2h) for the problem in **y**. The latter is almost trivial problem and of course former is easy to verify a linear program (LP). We know that if Lagrangian subproblems have so-called *integer polyhdron*[1] then the optimal dual value $z_D^*$ is equal to the optimal value if integer restriction was relaxed in the original model (Guignard, 2003, Cor. 7.1). Hence, on solving the dual problem, we only end up with a dual value equivalent to $z_{LP}^*$, which indicates no improvements in problem solving.

2. Decompose over $\mathcal{B}$: it requires to relax (2f) and (2d) due to sum over $j$. Furthermore, the constraint (2g) needs to be relaxed as well since the variable $y_{kt}$ links all $j \in \mathcal{B}$. Hence, we still have a huge dimension of the Lagrangian multiplier, and the resulting sub-problems can be decomposed over $j \in \mathcal{B}$ and variable **y**. The problem in **y** is clearly trivial where as the decomposed problem in $\mathbf{x}_{j,\cdot,\cdot}$ for each $j \in \mathcal{B}$ subject to constraints (2b)-(2c) and (2e) results in an integer polyhedron. A simple intuition is simply for each $j \in \mathcal{B}_\ell$, construct a graph with a source node of supply of $p_j$, a layer of nodes where each node corresponds to two-tuples $(j,t)$, where $t \in \mathcal{T}_\ell$. Finally, create another layer of nodes representing each $k \in \mathcal{N}$. Furthermore, a sink node is also needed. Arcs are allowed only from source to $(j,t)$ nodes and $(j,t)$ nodes to $k$ nodes and finally, all nodes corresponding to $k \in \mathcal{N}$ are linked to the sink node. The capacity limitation on arcs from $(j,t)$ to $k$ nodes is 1 and similarly, all incoming arcs to $(j,t)$ nodes has a capacity of one. The latter capacity limitation ensures constraint (2e) which is still in the model and the former capacity limitation is just because we have a binary variable and hence $x_{j,k,t} \leq 1$. Since it is reduced to a min-cost network flow problem we have $z_{LP}^* = z_D^*$.

---
[1] all extreme points of the corresponding polyhedron are integers

3. Lastly, we observe that decomposing over $k \in \mathcal{N}$ requires relaxing constraints (2c), (2b), (2e), and (2h) from the original model. The decomposition over $k \in \mathcal{N}$ leads to a subproblem that is nearly trivial (see Corollary 4.2.1)

## 4. Variable-layering-based Lagrangian dualization

In this section, we present a variable-layering-based Lagrangian dualization algorithm to address aforementioned challenges and solve the model (2) efficiently.

The basic idea is that subproblems should be established such that most of the important polyhedral properties are conserved. To that end, we must deal with constraints (2f)-(2g) with caution as they couple time periods as well as batteries. Hence, we use so-called *layering strategies* to create exploitable structure. The techniques of *variable layering* (Glover and Klingman, 1988), variable splitting (Jörnsten and Näsberg, 1986) and Lagrangian decomposition (Guignard, 2003) are highly similar to each other and have been successfully used in various decomposition techniques (Holmberg, 1998; Zhao et al., 2018; Granfeldt et al., 2023). The general idea is that one must create a copy of existing variables and add a so-called copy-constraint that can then be dualized using a Lagrangian multiplier. We explain this by implementing it for our model (2). The resultant reformulation will lead to two subproblems, one producing tighter bounds than standard relaxation techniques while the other having an integer polyderon, as will be proved in following propositions. Those properties of the meticulously designed subproblems brings superior computational performance, and thus needs be articulated.

Specifically, for a fixed value of $0 < \beta < 1$, we define an equivalent model to (2). Firstly, we define two polyhedral sets as follows:

$$P^{\mathbf{x}} := \left\{ x_{jkt} \in \{0,1\}, j \in \mathcal{B}_\ell, k \in \mathcal{N}, t \in \mathcal{T}_\ell, \ell \in \mathcal{L} \middle| \text{(2b)--(2e) satisfied} \right\}, \tag{3a}$$

$$P^{\mathbf{y}\text{-s}} := \left( (\mathbf{y}, \mathbf{s}) \middle| \begin{array}{l} y_{kt} \geq \left( \sum_{\ell \in \mathcal{L}} \sum_{j \in \mathcal{B}_\ell} s_{jkt+1} - \sum_{\ell \in \mathcal{L}} \sum_{j \in \mathcal{B}_\ell} s_{jkt} \right), \forall k, t, \\ y_{kt} \geq (s_{jkt} - s_{jkt+1}), \forall j, k, t, \ell, \\ \sum_{k \in \mathcal{N}} y_{kt} \leq \gamma, \forall t \in \mathcal{T} \setminus \{T\}, \\ y_{kt} \in \{0,1\}, \forall k \in \mathcal{N}, t \in \mathcal{T} \setminus \{T\}, \\ s_{jkt} \in \{0,1\}, \forall j \in \mathcal{B}_\ell, k \in \mathcal{N}, t \in \mathcal{T}_\ell, \ell \in \mathcal{L}. \end{array} \right) \tag{3b}$$

An equivalent model to (2) is stated as follows:

$$z_e^* := \min \left( \beta \sum_{\ell \in \mathcal{L}} \sum_{t \in \mathcal{T}_\ell} c_t \left( \sum_{j \in \mathcal{B}_\ell} \sum_{k \in \mathcal{N}} x_{jkt} \right) + w \sum_{k \in \mathcal{N}} \sum_{t \in \mathcal{T}} y_{kt} + (1-\beta) \sum_{\ell \in \mathcal{L}} \sum_{t \in \mathcal{T}_\ell} c_t \left( \sum_{j \in \mathcal{B}_\ell} \sum_{k \in \mathcal{N}} s_{jkt} \right) \right) \tag{4a}$$

s.t. $\mathbf{x} = \mathbf{s},$ (4b)

$\mathbf{x} \in P^{\mathbf{x}}, (\mathbf{y}, \mathbf{s}) \in P^{\mathbf{y}\text{-s}}.$ (4c)

We ensure that $z_e^* = z^*$ (refer to (2a)) due to the copy-constraint (4b). Next we dualize/relax (4b) using a Lagrangian multiplier $\boldsymbol{\mu}$. The Lagrangian function is defined as

$$L(\mathbf{x}, \mathbf{y}, \mathbf{s}; \boldsymbol{\mu}) := \left( \beta \sum_{t \in \mathcal{T}} c_t \left( \sum_{j \in \mathcal{B}} \sum_{k \in \mathcal{N}} x_{jkt} \right) + w \sum_{k \in \mathcal{N}} \sum_{t \in \mathcal{T}} y_{kt} \right.$$
$$\left. + (1-\beta) \sum_{t \in \mathcal{T}} c_t \left( \sum_{j \in \mathcal{B}} \sum_{k \in \mathcal{N}} s_{jkt} \right) + \sum_{t \in \mathcal{T}} \sum_{k \in \mathcal{N}} \sum_{j \in \mathcal{B}} \mu_{jkt}(x_{jkt} - s_{jkt}) \right).$$

Since $L(\mathbf{x}, \mathbf{y}, \mathbf{s}; \boldsymbol{\mu})$ is a linear function and both sets $P^{\mathbf{x}}, P^{\mathbf{y}\text{-s}}$ are closed (linear equalities/inequalities) and bounded (all involved variables have an upper bound of 1), by the classic Weierstrass theorem an optimal solution exists. Hence, we define the Lagrangian dual function as $h(\boldsymbol{\mu}) := \min_{\mathbf{x} \in P^{\mathbf{x}}, (\mathbf{y},\mathbf{s}) \in P^{\mathbf{y}\text{-}\mathbf{s}}} L(\mathbf{x}, \mathbf{y}, \mathbf{s}; \boldsymbol{\mu})$. Our first course of action is to separate into two Lagrangian subproblems, one for the polyhedron $P^{\mathbf{x}}$ and the other for $P^{\mathbf{y}\text{-}\mathbf{s}}$ as follows:

$$h(\boldsymbol{\mu}) := \min_{\mathbf{x} \in P^{\mathbf{x}}, (\mathbf{y},\mathbf{s}) \in P^{\mathbf{y}\text{-}\mathbf{s}}} L(\mathbf{x}, \mathbf{y}, \mathbf{s}; \boldsymbol{\mu}) := h_{\mathbf{x}}(\boldsymbol{\mu}) + h_{\mathbf{y},\mathbf{s}}(\boldsymbol{\mu}), \tag{5a}$$

$$\text{where } h_{\mathbf{x}}(\boldsymbol{\mu}) := \min_{\mathbf{x} \in P^{\mathbf{x}}} L^1(\mathbf{x}; \boldsymbol{\mu}) = \min_{\mathbf{x} \in P^{\mathbf{x}}} \sum_{\ell \in \mathcal{L}} \sum_{t \in \mathcal{T}_\ell} \sum_{j \in \mathcal{B}_\ell} \sum_{k \in \mathcal{N}} (\beta c_t + \mu_{jkt}) x_{jkt}, \tag{5b}$$

$$\text{and } h_{\mathbf{y},\mathbf{s}}(\boldsymbol{\mu}) := \min_{(\mathbf{y},\mathbf{s}) \in P^{\mathbf{y}\text{-}\mathbf{s}}} L^2(\mathbf{y}, \mathbf{s}; \boldsymbol{\mu})$$

$$:= \min_{(\mathbf{y},\mathbf{s}) \in P^{\mathbf{y}\text{-}\mathbf{s}}} \left( \sum_{\ell \in \mathcal{L}} \sum_{t \in \mathcal{T}_\ell} \sum_{j \in \mathcal{B}_\ell} \sum_{k \in \mathcal{N}} \underbrace{((1-\beta)c_t - \mu_{jkt})}_{\tilde{c}_{jkt}} s_{jkt} \right.$$

$$\left. + w \sum_{k \in \mathcal{N}} \sum_{t \in \mathcal{T} \setminus \{T\}} y_{kt} \right). \tag{5c}$$

**Observation 4.1.** *Consider a feasible set $\tilde{P}^{\mathbf{y}\text{-}\mathbf{s}} := P^{\mathbf{y}\text{-}\mathbf{s}} \cap \left\{ \mathbf{s} \mid \sum_{\ell \in \mathcal{L}} \sum_{j \in \mathcal{B}_\ell} s_{jkt} \leq 1, \forall k \in \mathcal{N}, \forall t \in \mathcal{T} \right\}$. Given that $\tilde{P}^{\mathbf{y}\text{-}\mathbf{s}} \subseteq P^{\mathbf{y}\text{-}\mathbf{s}}$, it follows that $h_{\mathbf{y},\mathbf{s}}(\boldsymbol{\mu}) \geq \min_{(\mathbf{y},\mathbf{s}) \in \tilde{P}^{\mathbf{y}\text{-}\mathbf{s}}} L^2(\mathbf{y}, \mathbf{s}; \boldsymbol{\mu}) =: \tilde{h}_{\mathbf{y},\mathbf{s}}(\boldsymbol{\mu})$. This inequality is likely to be strict. The only downside is the addition of $N \cdot T$ constraints to the $\mathbf{y}, \mathbf{s}$ Lagrangian subproblem. Our computational experiments have shown that the tighter set $\tilde{P}^{\mathbf{y}\text{-}\mathbf{s}}$ leads to faster convergence without significant computational burden. Therefore, we define $\tilde{h}(\boldsymbol{\mu}) := h_x(\boldsymbol{\mu}) + \tilde{h}_{\mathbf{y},\mathbf{s}}(\boldsymbol{\mu})$. Most importantly, this facilitates a specific reformulation to be derived later.*

We solve the dual problem as the Lagrangian dual function $\tilde{h}_D^* = \max_{\boldsymbol{\mu}} \tilde{h}(\boldsymbol{\mu})$ to obtain a lower bound to the optimal value of our model (2). Note that since we do not have strong-duality due to binary variables it is important to understand the quality of lower bound $\tilde{h}_D^*$ as compared to the lower bound obtained from linear relaxation of (2) i.e. $z_{\text{LP}}^*$ and possibly other classical Lagrangian relaxations. From proposition in (Jörnsten and Näsberg, 1986, p. 318) it is clear that the newly obtained lower bound from variable splitting/layering is also at least as strong as that for the traditional Lagrangian relaxation approaches (discussed in Sec. 3.1). However, to show that the last strict inequality holds in $z_e^* = z^* \geq \tilde{h}_D^* > z_{\text{LP}}^*$, we must prove that at least one of the two sub-problems corresponding to subproblems in $P^{\mathbf{x}}$ and $\tilde{P}^{\mathbf{y}\text{-}\mathbf{s}}$ do not have integer extreme points (Guignard, 2003, Cor. 7.1).

**Proposition 4.2** ($\tilde{P}^{\mathbf{y}\text{-}\mathbf{s}}$ *is not an integer polyhderon*). *The subproblem in $\tilde{P}^{\mathbf{y}\text{-}\mathbf{s}}$ can be polynomially reduced to a multiple knapsack problem (MKP). Hence, $\tilde{P}^{\mathbf{y}\text{-}\mathbf{s}}$ does not have integer extreme points.*

*Proof.* We provide a special instance of the aforementioned Lagrangian sub-problem and show that it is equivalent to any general multiple knapsack problem. Letting, for $k \in \{1, \ldots, n\}$ and $t \in \{1, \ldots, m\}$, $x_{kt} = 1$ if item $k$ is included in the knapsack $t$, and $x_{kt} = 0$ otherwise, the MKP can be modeled as:

$$z := \max \left\{ \sum_{k=1}^{n} \sum_{t=1}^{m} v_{kt} x_{kt} \, \middle| \, \sum_{k=1}^{n} x_{kt} \leq C_t, \forall t \in \{1, \ldots, m\}, x_{kt} \in \{0,1\}, \forall k \in \{1, \ldots, n\}, t \in \{1, \ldots, m\} \right\}, \tag{6}$$

where $v_{kt} \geq 0$ is a cost coefficient and $C_t > 0$ is the weight restriction on each knapsack $t$. Consider a special case of the Lagrangian sub-problem i.e. $\min_{(\mathbf{y},\mathbf{s}) \in \tilde{P}^{\mathbf{y}\text{-}\mathbf{s}}} L^2(\mathbf{y}, \mathbf{s}; \boldsymbol{\mu})$ when $w = 0$. We denote a new coefficient for each $k \in \mathcal{N}$, $\ell \in \mathcal{L}, t \in \mathcal{T}_\ell$ as

$$\alpha_{kt} := \begin{cases} 0, & \text{if } (1-\beta)c_t - \mu_{jkt} \geq 0 \text{ for all } j \in \mathcal{B}_\ell, \\ \min\{(1-\beta)c_t - \mu_{jkt} \mid j \in \mathcal{B}_\ell\}, & \text{otherwise.} \end{cases}$$

We denote $n := N$ (# of CBs), $m := T$ (# of time periods), $v_{kt} = -\alpha_{kt}$, and $C_t = \gamma$ then the Lagrangian sub-problem for $\mathbf{y}, \mathbf{s}$ can be solved by solving the following equivalent MKP [2]:

---
[2] note that having a homogeneous knapsacks is still $\mathcal{NP}$-hard

$$h := -\max\left\{\sum_{k=1}^{n}\sum_{t=1}^{m} v_{kt}y_{kt} \,\middle|\, \sum_{k=1}^{n} y_{kt} \leq C_t, \forall t = 1, \ldots, m\right\}. \quad (7)$$

The instance of the Lagrangian sub-problem in $\mathbf{y}$ as represented in (7), is a MKP, as expressed in (6). An optimal solution $(\mathbf{y}^*)$ to (7) can be used to deduce optimal $s^*_{jkt} = 1$, when $j \in \arg\min_{j \in \mathcal{B}_\ell}\{(1-\beta)c_t - \mu_{jkt}|\alpha_{kt} < 0\}$, otherwise $s^*_{jkt} = 0$. It is easy to follow that $(\mathbf{y}^*, \mathbf{s}^*) \in \arg\min\left\{L^2(\mathbf{y}, \mathbf{s}; \boldsymbol{\mu})|(\mathbf{y}, \mathbf{s}) \in \tilde{P}^{\mathbf{y}\text{-}\mathbf{s}}\right\}$. Since MKP is $\mathcal{NP}$-complete, the proposition follows. □

**Corollary 4.2.1.** *In Sec. 3.1 we discussed relaxing constraints* (2c), (2b), (2e), *and* (2h) *from the original model to decompose over* $k \in \mathcal{N}$. *This has the same constraints as* $\tilde{P}^{\mathbf{y}\text{-}\mathbf{s}}$ *except that constraint* (2h) *is also relaxed. Hence, the resulting subproblem is a LP and cannot be reduced to MKP. Hence, we get weak lower bound equivalent to relaxing binary restriction in the original model.*

Although we have already established that one of the subproblems is $\mathcal{NP}$-hard, we further demonstrate a crucial property of $P^{\mathbf{x}}$ that enhances the solution's effectiveness in computational tests.

**Proposition 4.3** ($P^{\mathbf{x}}$ *is an integer polyderon*). *The polyhedron $P^{\mathbf{x}}$ has integer extreme points and corresponding subproblem in* $\mathbf{x}$ *is equivalent to solving a min-cost network flow problem. (see appendix for proof)*

To summarize, at each iteration of solving the Lagrangian dual problem, we address two Lagrangian sub-problems: one in the variable space $\mathbf{x}$ and the other in the variable space $\mathbf{y}, \mathbf{s}$. As shown in Prop. 4.3, the subproblem in $\mathbf{x}$ can be reduced to a min-cost network flow problem and thus solved in polynomial time. However, the problem in $(\mathbf{y}, \mathbf{s})$ poses significant computational burden, as proved in Prop. 4.2, leading to longer computational time. Thus, the problem in $(\mathbf{y}, \mathbf{s})$ constitutes the computational bottleneck. To address this bottleneck, we propose a new formulation for the problem in $(\mathbf{y}, \mathbf{s})$ in upcoming section.

### 4.1. Polyhedral analysis of the sub-problem in $(\mathbf{y}, \mathbf{s})$ and time-block reformulations

In this subsection, we reformulate the Lagrangian subproblem in the variable space $(\mathbf{y}, \mathbf{s})$ to further improve computational performance. First, we introduce a discovered property of the optimal solution of the subproblem, which led to a new so-called *time block* reformulation. Let us consider an instance with $T = 5$ and $\mathcal{B} = \{1, 2, 3, 4\}$ and create a network as illustrated Fig. 2 for any given $k \in \mathcal{N}$, and assume we know $\bar{\mathbf{y}} \in \{y_{kt} \in \{0, 1\}, k \in \mathcal{N}, t \in \mathcal{T} \setminus \{T\} \mid \sum_{k \in \mathcal{N}} y_{kt} \leq \gamma, t \in \mathcal{T} \setminus \{T\}\}$. Each round nodes correspond to $(j, t), \ell \in \mathcal{L}, j \in \mathcal{B}_\ell, t \in \{1, \ldots, T\} =: \mathcal{T}$ and the number on top of the node (in red) is $\bar{c}_{jkt} := (1-\beta)c_t - \mu_{jkt}$, where $\mu_{jkt}$ is the Lagrangian multiplier at a given iteration. The square nodes refer to time period $t \in \{0\} \cup \mathcal{T}$ when $\bar{y}_{kt} = 1$ or $t = 0$ or $t = T$ when all charging activities end. Note that $t = 0$ implies start of first time period and there is no limitation on the number of switchings at the beginning. Same can be said about the last time period $t = T$.

For any given $\bar{\mathbf{y}}$, we can draw a similar graph as in the Figure 2. $\sum_{\ell \in \mathcal{L}} \sum_{j \in \mathcal{B}_\ell} s_{jkt} \leq 1, k \in \mathcal{N}, t \in \mathcal{T}$ implies that for time periods $\{1, 2, 3\}$, i.e. time block $(0, 3)$, we can select only one battery $j \in \mathcal{B}_\ell, \ell \in \mathcal{L}$ to charge in battery bay . Since we aim for reducing the cost, i.e., sum of red numbers in each row of a block in the figure, it is obvious that we should select $j = 4$ for the first time block $(0, 3)$ and $j = 3$ for the time block $(3, 5)$ for a given $k \in \mathcal{N}$. In Prop. 4.4, we formally present the characterization of optimal solution of the Lagrangian subproblem in $(\mathbf{y}, \mathbf{s})$ used to deduce the new formulation.

**Proposition 4.4** (*Properties of optimal solution to the Lagrangian sub-problem in* $(\mathbf{y}, \mathbf{s})$). *Let* $(\mathbf{y}^*, \mathbf{s}^*)$ *be an optimal solution of* $\min_{(\mathbf{y}, \mathbf{s}) \in \tilde{P}^{\mathbf{y}\text{-}\mathbf{s}}} L^2(\mathbf{y}, \mathbf{s}; \boldsymbol{\mu})$. *For a given pair* $t_1, t_2$ *and* $k \in \mathcal{N}$, *where* $t_1 \in \{0, \ldots, T-1\}$ *and* $t_2 \in \{1, \ldots, T\}$ *such that* $t_2 > t_1$ *and* $y^*_{kt_2} = y^*_{kt_1} = 1$ *and* $\sum_{t=t_1+1}^{t_2-1} y^*_{kt} = 0$ *then following holds*[3]:

- *For any* $j \in \mathcal{B}_\ell, \ell \in \mathcal{L}$ *and* $t_2 \leq n_\ell$ *(i.e. in the right demand window) if* $\sum_{t=t_1+1}^{t_2}((1-\beta)c_t - \mu_{jkt}) \geq 0$ *then following holds* $s^*_{jkt_1+1} = s^*_{jkt_1+2} = \ldots = s^*_{jkt_2} = 0$

---
[3]Assume initial condition $y^*_{k0} = y^*_{kT} = 1$

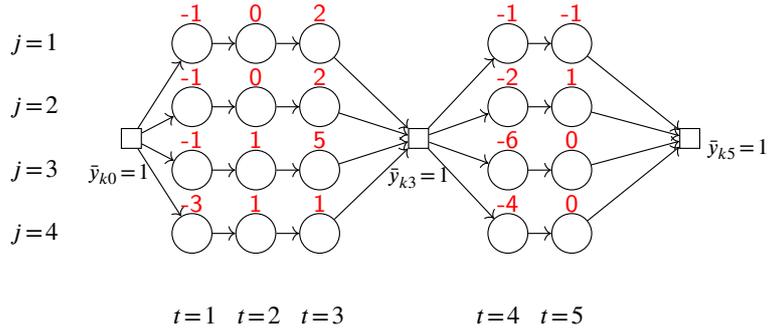

**Figure 2:** An example of the time block representation of the Lagrangian sub-problem in $(\mathbf{y}, \mathbf{s})$.

- *For any $j^* \in \mathcal{B}_\ell, \ell \in \mathcal{L}$ and $t_2 \leq n_\ell$ if $j^* \in \arg\min_{j \in \mathcal{B}} \sum_{t=t_1+1}^{t_2}((1-\beta)c_t - \mu_{jkt})$ such that $((1-\beta)c_t - \mu_{j^*kt}) < 0$ then following holds: $s^*_{j^*kt_1+1} = \ldots = s^*_{j^*kt_2} = 1$.*

*Proof.* For a fixed value of $\mathbf{y}^*$ (that satisfies (2h) and (2j)), the Lagrangian subproblem in $(\mathbf{y}, \mathbf{s})$ is reduced to following:

$$\min \sum_{\ell \in \mathcal{L}} \sum_{t \in \mathcal{T}_\ell} \sum_{j \in \mathcal{B}_\ell} \sum_{k \in \mathcal{N}} \bar{c}_{jkt} s_{jkt}, \tag{8a}$$

$$y^*_{kt} \geq \left( \sum_{\ell \in \mathcal{L}} \sum_{j \in \mathcal{B}_\ell} s_{jkt+1} - \sum_{\ell \in \mathcal{L}} \sum_{j \in \mathcal{B}_\ell} s_{jkt} \right), \quad k \in \mathcal{N}, t \in \mathcal{T} \setminus \{T\}, \tag{8b}$$

$$y^*_{kt} \geq (s_{jkt} - s_{jkt+1}), \quad \forall j \in \mathcal{B}_\ell, k \in \mathcal{N}, t \in \mathcal{T}_\ell \setminus \{T\}, \ell \in \mathcal{L}, \tag{8c}$$

$$\sum_{\ell \in \mathcal{L}} \sum_{j \in \mathcal{B}_\ell} s_{jkt} \leq 1, \quad \forall k \in \mathcal{N}, t \in \mathcal{T}, \tag{8d}$$

$$s_{jkt} \in \{0, 1\}, \quad \forall j \in \mathcal{B}_\ell, k \in \mathcal{N}, t \in \mathcal{T}_\ell, \ell \in \mathcal{L}, \tag{8e}$$

where $\bar{c}_{jkt} := (1-\beta)c_t - \mu_{jkt}$. We note that if for a given $k, t$ we have $y^*_{kt} = 0$, it implies $s_{jkt} \leq s_{jkt+1}$ (8c) and if we sum both sides we get $\sum_{\ell \in \mathcal{L}} \sum_{j \in \mathcal{B}_\ell} s_{jkt} \leq \sum_{\ell \in \mathcal{L}} \sum_{j \in \mathcal{B}_\ell} s_{jkt+1}$. Due to previous statement and (8b) we get that $s_{jkt+1} = s_{jkt}, \forall j \in \mathcal{B}$ if $y^*_{kt} = 0$. Similarly, if $y^*_{kt} = 1$, the constraint (8c) is redundant. Furthermore, (8d) makes (8b) redundant. So, given a $\mathbf{y}^*$ satisfying (2h) and (2j) we can split the time horizon for each $k \in \mathcal{N}$ into time blocks/intervals as follows

$$L^k(\mathbf{y}^*) := \{t \mid y^*_{kt} = 1, t \in \mathcal{T} \setminus \{T\}\} \cup \{0\} \cup \{T\}, k \in \mathcal{N} \tag{9a}$$

$$Z^k(\mathbf{y}^*) := \left\{ (t_i, t_{i+1}) \mid t_i, t_{i+1} \in L^k(\mathbf{y}^*) \text{ such that } t_i < t_{i+1}; \nexists t \in L^k(\mathbf{y}^*) \text{ such that } t_i < t < t_{i+1} \right\}, k \in \mathcal{N} \tag{9b}$$

$$Z^k_<(\mathbf{y}^*) := \left\{ (t_1, t_2) \in Z^k(\mathbf{y}^*) \,\Big|\, \min_{j \in \mathcal{B}} \left\{ \sum_{t=t_1+1}^{t_2} \bar{c}_{jkt} \right\} < 0 \right\}, k \in \mathcal{N} \tag{9c}$$

For each $(t_1, t_2)$ in $Z^k_<(\mathbf{y}^*)$ we can only select one of the batteries $j \in \mathcal{B}$. We compute $j^* \in \arg\min_{j \in \mathcal{B}} \left\{ \sum_{t=t_1+1}^{t_2} \bar{c}_{jkt} \right\}$ and set $s^*_{j^*kt_1+1} = \ldots = s^*_{j^*kt_2} = 1$ and $s^*_{jkt_1+1} = \ldots = s^*_{jkt_2} = 0, j \in \mathcal{B} \setminus \{j^*\}$. The proposition follows. □

Using results from Prop. 4.4, we propose the *time-block formulation*, which aims to enhance computational efficiency through discrete time segmentation. We define the set of all possible time blocks for each port $k \in \mathcal{N}$:

$$\mathcal{T}^k := \{(t_i, t_{i+1}) \mid t_i \in \{0\} \cup \mathcal{T} \setminus \{T\}, t_{i+1} \in \mathcal{T}, t_i < t_{i+1}\}, \quad k \in \mathcal{N}, \tag{10}$$

A subset $\mathcal{T}_<^k$ comprises time blocks where the respective objective function's coefficient is $\bar{c}_{jkt} = (1-\beta)c_t - \mu_{jkt}$ is negative:

$$\mathcal{T}_<^k := \left\{ (t_1, t_2) \in \mathcal{T}^k \;\Big|\; \min_{j \in \mathcal{B}_\ell, \ell \in \mathcal{L}} \left\{ \sum_{t=t_1+1}^{t_2} \bar{c}_{jkt} \right\} < 0 \right\}, \quad k \in \mathcal{N}, \tag{11}$$

with $\bar{c}_{jkt} := (1-\beta)c_t - \mu_{jkt}$ representing the adjusted cost parameters. Next we state the new formulation for solving the Lagrangian subproblem in $(\mathbf{y}, \mathbf{s})$ by creating a new equivalent model in variable space $\mathbf{y}, \lambda$. Here, $\lambda_{k,t_1,t_2}$ equals one if a port $k \in \mathcal{N}$ charges a battery $j \in \mathcal{B}$ such that $j \in \arg\min_{j \in \mathcal{B}_\ell, \ell \in \mathcal{L}} \left\{ \sum_{t=t_1+1}^{t_2} \bar{c}_{jkt} \mid (t_1, t_2) \in \mathcal{T}_<^k \right\}$ and otherwise zero. Since there is no index corresponding to $j \in \mathcal{B}$ in model (12) we pre-compute the cost of using a time-block using the coefficient $d_{t_1 t_2}^k := \min_{j \in \mathcal{B}_\ell, \ell \in \mathcal{L}} \left\{ \sum_{t=t_1+1}^{t_2} \bar{c}_{jkt} \mid (t_1, t_2) \in \mathcal{T}_<^k \right\}$ used in the objective function as well (12a). The variables $\mathbf{y}$ represent switchings as used in the original model (2).

Constraints (12b) prevent overlapping time blocks for $k \in \mathcal{N}$. For instance, for any given $k$ if $y_{k3} = y_{k6} = y_{k10} = 1$ then if time blocks considered in $\mathcal{T}^k$ are $\lambda_{k,0,3} = \lambda_{k,3,6} = \lambda_{k,6,10} = \lambda_{k,10,24} = 1$ (Fig. 3). The constraint (12b) simply imposes a restriction that if $\lambda_{k,0,3} = 1$ then $\lambda_{k,1,2} = 0$ and similarly if $\lambda_{k,3,6} = 1$ then $\lambda_{k,4,5} = 0$. This ensures we do not select more than one $j$ for the same time period which is a constraint in $\tilde{P}^{\mathbf{y}\text{-}\mathbf{s}}$.

Constraints (12c) and (12d) link the start and end of time blocks to variables $\mathbf{y}$. For instance, if a time block is selected with starting time period $t_1 \in \mathcal{T} \setminus \{T\}$ then $y_{k,t_1} \geq 1$ is enforced, hence, $y_{kt_1}$ represents switching. Another utility of these constraints is that since $y_{k,t_1} \leq 1$ the constraint also implies a restriction that starting time period $t_1$ only one time block can have $t_1$ as the starting point. Constraint (12e) limits the number of switchings at any time to $\gamma$.

The resulting optimization problem is reformulated as (Time-block formulation $F_{\text{tb}}$):

$$\min \sum_{k \in \mathcal{N}} \sum_{(t_1, t_2) \in \mathcal{T}_<^k} d_{t_1 t_2}^k \lambda_{k t_1 t_2} + w \sum_{k \in \mathcal{N}} \sum_{t \in \mathcal{T} \setminus \{T\}} y_{kt}, \tag{12a}$$

$$\text{s.t.} \quad \lambda_{k t_1 t_2} + \lambda_{k t_1' t_2'} \leq 1, \qquad \forall k \in \mathcal{N}, (t_1, t_2), (t_1', t_2') \in \mathcal{T}_<^k, \; t_1' < t_2 \wedge t_2' > t_1, \tag{12b}$$

$$y_{k t_1} \geq \sum_{t_1+1}^{T} \lambda_{k t_1 t}, \qquad \forall k \in \mathcal{N}, t_1 \in \{1, \ldots, T-1\}, \tag{12c}$$

$$y_{k t_2} \geq \sum_{t=0}^{t_2-1} \lambda_{k t t_2}, \qquad \forall k \in \mathcal{N}, t_2 \in \{1, \ldots, T-1\}, \tag{12d}$$

$$\sum_{k \in \mathcal{N}} y_{kt} \leq \gamma, \qquad \forall t \in \mathcal{T} \setminus \{T\}, \tag{12e}$$

$$\lambda_{k t_1 t_2} \in \{0, 1\}, \qquad \forall k \in \mathcal{N}, t \in \mathcal{T}, (t_1, t_2) \in \mathcal{T}_<^k, \tag{12f}$$

$$y_{kt} \in \{0, 1\}, \qquad k \in \mathcal{N}, t \in \mathcal{T} \setminus \{T\}. \tag{12g}$$

Now let us call the feasible set defined by constraints in the time-block formulation (model (12)) as $\tilde{P}^{\mathbf{y}\text{-}\lambda}$. There are several crucial points that merit discussion. Firstly, it is essential to evaluate whether the new time-block formulation for a given $\mu$ is equivalent to the orginal subproblem in $(\mathbf{y}, \mathbf{s})$ despite having a different variable space $(\mathbf{y}, \lambda)$. Secondly, when proposing a new formulation, discussing its compactness and tightness becomes vital. Thirdly, theoretical claims should be substantiated with experimental evidence from realistic instances of interest. The proofs of equivalence (Prop. 4.5) and tightness (Prop. 13) are technical and hence, described in detail in the appendix. An overview is provided in Fig. 4.

**Proposition 4.5** (Equivalence of $F_r$ and $F_{\text{tb}}$). *The time block formulation $F_{\text{tb}}$ admits an optimal solution $(\mathbf{y}^*, \lambda^*) \in \tilde{P}^{\mathbf{y}\text{-}\lambda}$. This solution can be correspondingly mapped to $(\mathbf{y}^*, \mathbf{s}^*) \in \tilde{P}^{\mathbf{y}\text{-}\mathbf{s}}$ and retains optimality under the formulation $F_r$. (proof in the appendix)*

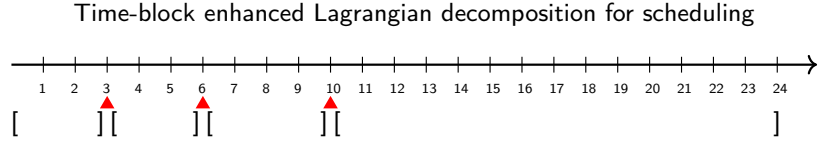

**Figure 3:** An example of disjunctive time-blocks for a given port $k = 1$ and known $\bar{y}_{1t} = 1, t = \{3, 6, 10\}$

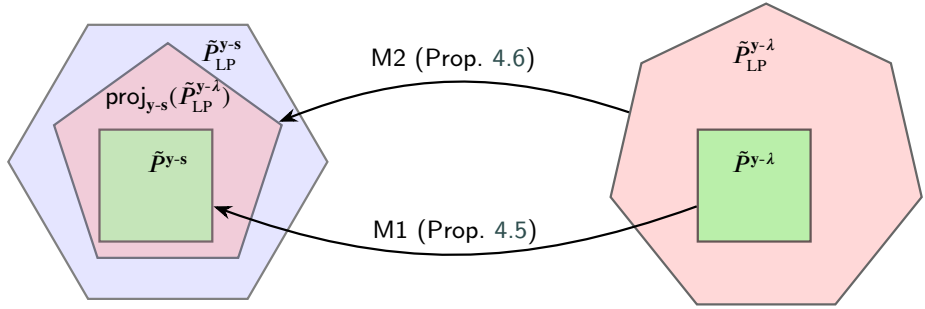

**Figure 4:** An illustration of mapping of the two formulations

Next, we discuss whether the new formulation $F_{\text{tb}}$ is more compact than $F_{\text{r}}$. Computing the dimension of $\lambda$ is challenging as it depends on the set $\mathcal{T}_<^k$, which varies with the Lagrangian multipliers and changes in each iteration. However, an upper bound for the dimension can instead be based on $\mathcal{T}^k$ (10). Consequently, $|\mathcal{T}_<^k| \leq \frac{T(T-1)}{2}$, a calculation derived from the sum of an arithmetic series. Therefore, the dimension of $\lambda$ is less than or equal to $N \cdot \frac{T(T-1)}{2}$. This should now be compared to the dimension of $\mathbf{s}$, which is $B \cdot N \cdot T$. Given that $T < B$ in all instances (in fact $T$ remains the same for all the instances), we have effectively reduced the number of variables in the new formulation $F_{\text{tb}}$. Notably, the set $\mathcal{T}_<^k$ is much smaller than $\mathcal{T}^k$, primarily because the value of $\mu$ is unlikely to be high enough to make all components of $\bar{c}$ negative over many iterations, as dictated by the mechanics of subgradient optimization.

Due to constraint (12b) it is again hard to compare the number of constraints of the two formulations. Note that the number of constraints corresponding (12b) grows exponentially as a function of the cardinality of $\mathcal{T}_<^k$ so it is not fair to do a worst case analysis for these sets of constraints. Finally, we compare the tightness of the two formulations which hugely impacts the computation time when using any branch and bound and cut solvers like Gurobi (Bixby et al., 2000). Let us consider the two linear polyhedra provided to us in the two formulations as follows:

$$\tilde{P}_{\text{LP}}^{\mathbf{y}\text{-}\lambda} := \left\{ (\mathbf{y}, \lambda) \mid \text{(12b)-(12e)}, \lambda_{jkt_1t_2} \in [0,1], y_{kt} \in [0,1] \right\}, \quad \tilde{P}_{\text{LP}}^{\mathbf{y}\text{-}\mathbf{s}} := \left\{ (\mathbf{y}, \mathbf{s}) \mid \text{(8b)-(8d)}, \sum_{k \in \mathcal{N}} y_{kt} \leq \gamma, t \right\} \tag{13a}$$

In order compare the new formulation to the old formulation we need to project $\tilde{P}_{\text{LP}}^{\mathbf{y}\text{-}\lambda}$ onto variable space $(\mathbf{y}, \mathbf{s})$. The projection can be defined as:

$$\text{proj}_{\mathbf{y},\mathbf{s}}(\tilde{P}^{\mathbf{y}\text{-}\lambda}) := \{(\mathbf{y}, \mathbf{s}) \mid \mathbf{s} = \text{M1}(\lambda), (\mathbf{y}, \lambda) \text{ satisfy (12b)-(12e)}\}, \tag{14}$$

where M1 is a mapping of $\lambda$ to $\mathbf{s}$. Hence, we need to show that $\text{proj}_{\mathbf{y},\mathbf{s}}(\tilde{P}_{\text{LP}}^{\mathbf{y}\text{-}\lambda}) \subset \tilde{P}_{\text{LP}}^{\mathbf{y}\text{-}\mathbf{s}}$. The strict inequality is important to ensure that we obtain a tighter formulation. For more background in comparing formulations, we recommend (Conforti et al., 2014, Ch. 2.2). Following deals with tightness of the new formulation:

**Proposition 4.6.** *Show that $\text{proj}_{\mathbf{y},\mathbf{s}}(\tilde{P}_{\text{LP}}^{\mathbf{y}\text{-}\lambda}) \subset \tilde{P}_{\text{LP}}^{\mathbf{y}\text{-}\mathbf{s}}$. (proof in the appendix)*

With the above propositions, we proved that the new formulation $F_{\text{tb}}$ use fewer variables and is tighter than the original formulation $F_{\text{r}}$. Now, we present practical evidence of the benefits obtained. Figures 5a and 5c illustrate the lower bounds and computation times for Lagrangian subproblems using the old formulation $F_{\text{r}}$ with instances having 100 and 300 batteries, respectively. Conversely, Figures 5b and 5d display the lower bounds and computation times for solving subproblems under the new formulation $F_{\text{tb}}$, also for instances with 100 and 300 batteries.

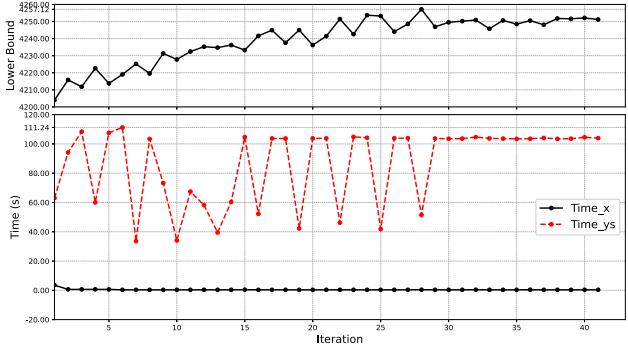

(a) Original formulation $F_r$ with $B = 100$

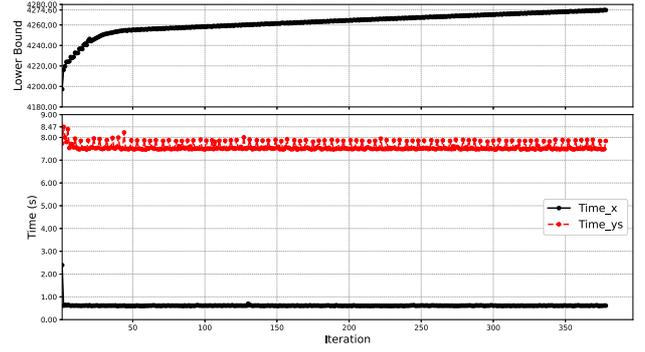

(b) Time-block formulation $F_{tb}$ with $B = 100$

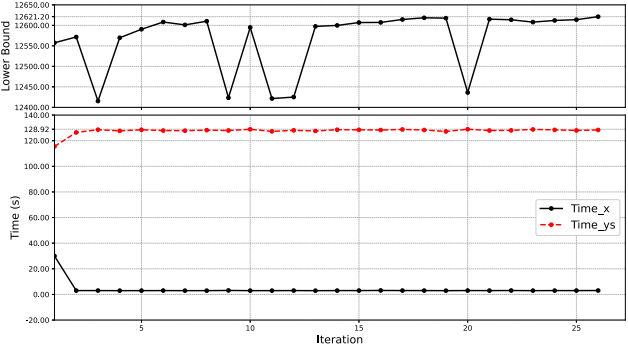

(c) Original formulation with $B = 300$

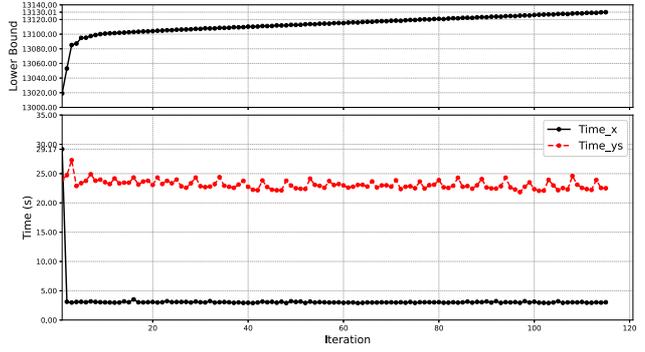

(d) Time-block formulation with $B = 300$

**Figure 5:** Impact of time-block formulation $F_{tb}$ on lower bounds and solution time to solve the subproblem in $(\mathbf{y}, \mathbf{s})$.

Notably, the computation time for solving subproblems in $\mathbf{x}$ is negligible[4], as expected (Prop. 4.3). Furthermore, the new time-block formulation achieves significantly higher lower bounds in both instances. This improvement is attributed to the faster solving times for subproblems in $\mathbf{y}, \mathbf{s}$ using the new formulation, allowing for more iterations within the same computational timeframe, thereby yielding better lower bounds. Note that a computational time limit of 100 seconds is imposed on each subproblem during an iteration. The total time limit is 3600 seconds. Consequently, the flat line observed for instance with 300 batteries (in Fig. 5c) around 120 seconds also accounts for additional model load time. This suggests that the model did not always converge to the optimal solution. Note that this reformulation impacts both the lower and upper bounds. As we have not yet introduced upper bounding techniques (i.e. primal feasibility menthods), the combined results are reserved for discussion later in Fig. 12.

## 5. Recovering feasible solutions and improving lower and upper bounds

To summarize, we have so far discussed the initial model defined for D-CSP (2), and then introduced variable layering to facilitate a Lagrangian decomposition into two independent subproblems (see (4), (5)). We have reformulated one of the subproblems to accelerate the solution process. Next, we will present an algorithm to obtain a feasible solution for the D-CSP model and will additionally apply local search techniques to enhance these feasible solutions. Furthermore, we will discuss the update of the Lagrangian multiplier. This will collectively form the components of subgradient optimization (SO) framework.

Firstly, we provide a brief overview of the framework summarized in Fig. 6. After relaxing the copy constraint (4b) in model (4) using an initial multiplier $\boldsymbol{\mu}^0$. The optimal solutions of the two subproblems are used by our variable fixing heuristic to compute feasible solutions and we also keep track so-called egodic iterates, which are essentially convex combinations of subproblem solutions (cf. Sec. 5.2). These ergodic iterates are used in the local search procedure.

---

[4]note that in the first iteration we implement a min-cost network using Google OR tools and hence, the first iteration has significantly higher computation time. In the later iterations we just update the edge weight and solve.

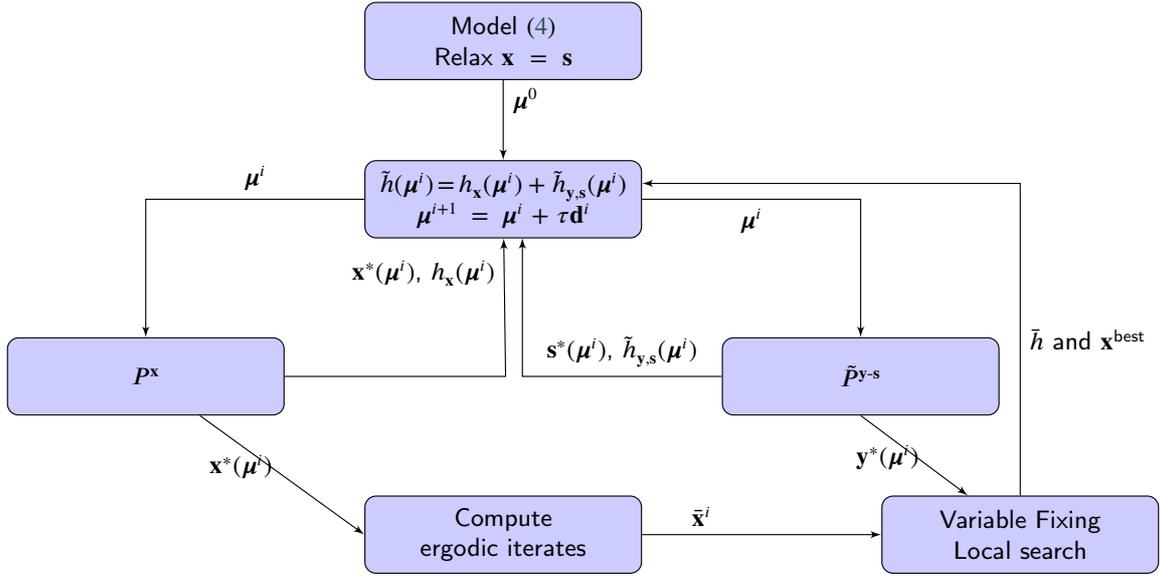

**Figure 6:** An overview of the Subgradient optimization procedure with variable fixing and local search.

After each iteration $i$, the new multiplier $\mu^{i+1}_{jkt}$ is computed as $\mu^{i+1}_{jkt} = \mu^i_{jkt} + \tau d^i_{jkt}$, where $\tau$—the step length according to Polyak's rule—is defined as $\tau = \theta^i \frac{\bar{h} - \tilde{h}(\mu^i)}{\|\mathbf{d}^i\|_2^2}$ with $\theta^i \in (0, 2)$. This adjustment is based on the difference between a known upper bound $\bar{h}$ and the dual value $\tilde{h}(\mu^i)$. The direction vector $\mathbf{d}^i$ can either be set to be equal to the subgradient $\mathbf{g}^i = (\mathbf{x}^*(\mu^i) - \mathbf{s}^*(\mu^i))$ or computed using MDS (1). The stopping criteria are either achieving $\epsilon$-optimality (gap)[5] with $\epsilon = 10^{-3}$ or hitting the 3600-second time limit, which is common in this study. Now, we present our novel variable fixing technique that takes in optimal subproblem solution for $\tilde{P}^{\mathbf{y}\text{-}\mathbf{s}}$ and computes a feasible solution.

### 5.1. A partial variable-fixing heuristic using disjunctive time-block-based bin-packing formulation

To compute a primal feasible solution, we utilize the optimal solution derived from solving the Lagrangian subproblem during each iteration. Given a subproblem solution $(\bar{\mathbf{y}}, \bar{\mathbf{s}}) \in \tilde{P}^{\mathbf{y}\text{-}\mathbf{s}}$ we consider the following optimization problem:

$$z^*(\bar{\mathbf{y}}) := \min \left( \sum_{\ell \in \mathcal{L}} \sum_{t \in \mathcal{T}_\ell} c_t \left( \sum_{j \in \mathcal{B}_\ell} \sum_{k \in \mathcal{N}} x_{jkt} \right) \,\middle|\, \mathbf{x} \in \mathrm{proj}_{P^{\mathbf{x}} | \bar{\mathbf{y}}}(A) \right),$$

where $A := \{(\mathbf{x}, \mathbf{y}) | (\mathbf{x}, \mathbf{y}) \text{ satisfies (2b)–(2j)}\}$ and for a given subproblem solution $\bar{\mathbf{y}}$, we define $\mathrm{proj}_{P^{\mathbf{x}} | \bar{\mathbf{y}}}(A) := \{\mathbf{x} \in P^{\mathbf{x}} | (\mathbf{x}, \bar{\mathbf{y}}) \in A\}$ is the projection onto variable space $\mathbf{x}$. Hence, in simple terms we fix $\mathbf{y} = \bar{\mathbf{y}}$ in model(2) and solve for $\mathbf{x}$. However, based on our computational experiments (to be presented later in Sec. 6), this approach proves to be inefficient and even results in no feasible solution within time limit. Therefore, we propose an alternative compact formulation to reduce computational time and find feasible solutions. Given an optimal subproblem solution $\bar{\mathbf{y}} = \mathbf{y}^*(\mu^i)$, our task is to find $\mathbf{x} \in P^{\mathbf{x}}$ such that $(\mathbf{x}, \bar{\mathbf{y}})$ is feasible in the original model(2)[6]. First, we split the time horizon for each port $k \in \mathcal{N}$ into disjunctive time-blocks that cover the entire horizon. For instance, for a given $k = 1$, we have $\bar{y}_{kt} = 0, t \in \{1, 2, 4, 5, 7, 8, 9, 11, 12, 13, 14, 15, 16, 17, 18, 19, 20, 21, 22, 23\}$ and the remaining $\bar{y}_{kt} = 1, t \in \{3, 6, 10\}$. This implies that at port $k = 1$ switch can take place at end of time periods 3, 6, and 10. These disjunctive time-blocks are illustrated in Fig. 3. and formally defined as

$$Z^k(\bar{\mathbf{y}}) := \left\{ (t_i, t_{i+1}) \,|\, t_i, t_{i+1} \in L^k(\bar{\mathbf{y}}) \text{ such that } t_i < t_{i+1}; \nexists t \in L^k(\bar{\mathbf{y}}) \text{ such that } t_i < t < t_{i+1} \right\}, k \in \mathcal{N}, \quad (15)$$

---

[5]optimality gap = $\frac{\bar{h}^{\text{best}} - \underline{h}^{\text{best}}}{\bar{h}^{\text{best}}}$, where $\bar{h}^{\text{best}}$ is the best-known upper bound and $\underline{h}^{\text{best}}$ is the best-known lower bound.

[6]it is possible that in some iterations the subproblem solution may not have any feasible $\mathbf{x}$.

**Algorithm 1** Subgradient Optimization with Variable Layering and Reformulations (SOVLR)

1: Initialize $\theta^0$, best estimate of objective value $\bar{h}^{\text{best}} = +\infty$ and lower bound $\underline{h}^{\text{best}} = -\infty$, heu = 20
2: $\boldsymbol{\mu}^0 =$ using known regression parameters $\omega_1^\ell, \omega_2^\ell$  ▷ See Sec. 5.3
3: **while** termination conditions not met **do**
4:     $\mathbf{x}^*(\boldsymbol{\mu}^i), \mathbf{y}^*(\boldsymbol{\mu}^i), \mathbf{s}^*(\boldsymbol{\mu}^i) =$ Solve Lagrangian sub-problems for $h_{\mathbf{x}}(\boldsymbol{\mu}^i)$ and $\tilde{h}_{(\mathbf{y},\mathbf{s})}(\boldsymbol{\mu}^i)$ in **parallel**  ▷ See Sec. 4.1
5:     **if** *Loc* & $\mathbf{x}^{\text{best}}$ is not None & $i \mod \frac{\text{heu}}{2} = 0$ **then**
6:         Use $\sigma, \tilde{\mathbf{x}}^{i-1}, \mathbf{y}^{\text{best}}$ and $\mathbf{x}^{\text{best}}$ to solve model(24) in **parallel** and update $\mathbf{x}^{\text{best}}, \underline{h}^{\text{best}}$ and $\mathbf{y}^{\text{best}}$
7:     **end if**
8:     $\underline{h}^{\text{best}} = \max\{h_{\mathbf{x}}(\boldsymbol{\mu}^i) + \tilde{h}_{\mathbf{y},\mathbf{s}}(\boldsymbol{\mu}^i), \underline{h}^{\text{best}}\}$
9:     Compute current subgradient $g_{jkt}^i := (x_{jkt}^i - s_{jkt}^i), \forall j, k, t$.
10:     **if** using *Simple* **then**
11:         Set $\mathbf{d}^i = \mathbf{g}^i$.
12:     **else**
13:         Update $\mathbf{d}^i$ using MDS method: $\mathbf{d}^i := \mathbf{g}^i + Y_{MDS}^i \mathbf{d}^{i-1}$.  ▷ See eq. (1)
14:     **end if**
15:     **if** $i \mod \text{heu} = 0$ **then**
16:         Update $\bar{h} =$ var-fix-binP$(\mathbf{y}^i)$  ▷ See model(16)
17:         $\bar{h}^{\text{best}} = \min\{\bar{h}, \bar{h}^{\text{best}}\}$
18:         Update $\mathbf{y}^{\text{best}}, \mathbf{x}^{\text{best}}$
19:     **end if**
20:     Compute step size $\tau$ (Polyak step length).
21:     **if** *Reg* **then**
22:         Compute regression parameters $\omega_1^\ell(\boldsymbol{\mu}^{i+\frac{1}{2}}), \omega_2^\ell(\boldsymbol{\mu}^{i+\frac{1}{2}})$ as indicated in eq. (18a)
23:         Update $\boldsymbol{\mu}^{i+1}$ as indicated in eq. (18b)
24:     **else**
25:         Update $\boldsymbol{\mu}^{i+1} = \boldsymbol{\mu}^i + \tau \mathbf{d}^i$.
26:     **end if**
27:     **if** *Loc* **then**
28:         Using previous ergodic iterate $\tilde{\mathbf{x}}^{i-1}$ compute $\tilde{\mathbf{x}}^i$ as shown in eq. (25)
29:     **end if**
30:     $i = i + 1$.
31: **end while**
32: Output the best estimate lower bound ($\underline{h}^{\text{best}}$) and best upper bound ($\bar{h}^{\text{best}}$).

where $L^k(\bar{\mathbf{y}}) := \{t \mid \bar{y}_{kt} = 1, t \in \mathcal{T} \setminus \{T\}\} \cup \{0\} \cup \{T\}$. For the example in Fig. 3 we have $Z^1(\bar{\mathbf{y}}) = \{(0,3),(3,6),(6,10),(10,24)\}$. For a subproblem optimal solution $\bar{\mathbf{y}}$, we present a new compact model

$$\min \sum_{k \in \mathcal{N}} \sum_{(t_1,t_2) \in Z^k(\bar{\mathbf{y}})} R_{t_1 t_2} \left\{ \sum_{\ell \in \mathcal{L}} \sum_{j \in \{\mathcal{B}_\ell | t_2 \leq n_\ell\}} \kappa_{t_1 t_2 j}^k \right\}, \tag{16a}$$

$$\text{s.t.} \sum_{\ell \in \mathcal{L}} \sum_{j \in \mathcal{B}_\ell | t_2 \leq n_\ell} \kappa_{t_1 t_2 j}^k \leq 1, \quad k \in \mathcal{N}, (t_1, t_2) \in Z^k(\bar{\mathbf{y}}), \tag{16b}$$

$$\sum_{k \in \mathcal{N}} \sum_{(t_1,t_2) \in Z^k(\bar{\mathbf{y}}) | t_2 \leq n_\ell, t \in [t_1+1, t_2]} \kappa_{t_1, t_2, j}^k \leq 1, \quad t \in \mathcal{T}, j \in \mathcal{B}_\ell, \ell \in \mathcal{L}, \tag{16c}$$

$$\sum_{k \in \mathcal{N}} \sum_{(t_1,t_2) \in Z^k(\bar{\mathbf{y}})} (t_2 - t_1) \kappa_{t_1, t_2, j}^k = p_j, \quad j \in \mathcal{B}_\ell, \ell \in \mathcal{L} \setminus \{L\}, \tag{16d}$$

$$\sum_{k \in \mathcal{N}} \sum_{(t_1,t_2) \in Z^k(\bar{\mathbf{y}})} (t_2 - t_1) \kappa_{t_1, t_2, j}^k = \lceil \alpha p_j \rceil, \quad j \in \mathcal{B}_L, \tag{16e}$$

$$\kappa_{t_1 t_2 j}^k \in \{0, 1\}, \quad (t_1, t_2) \in Z^k(\bar{\mathbf{y}}), k \in \mathcal{N}, j \in \mathcal{B}_\ell, t_2 \leq n_\ell, \ell, \tag{16f}$$

where $R_{t_1 t_2} := \sum_{t=t_1+1}^{t_2} c_t$ denotes the electricity cost associated with the time periods $t_1 + 1, \ldots, t_2$. The variable $\kappa_{t_1, t_2 j}^k$ equals 1 if time periods $t_1 + 1, \ldots, t_2$ on port $k$ is allocated to battery $j$. The constraint (16b) ensures that each port $k \in \mathcal{N}$ in each time block $(t_1, t_2) \in Z^k(\bar{\mathbf{y}})$ can be allocated to at most one battery. For instance, Fig. 7 illustrates a simplified example. Let us consider that we have the same set $Z^1(\bar{\mathbf{y}})$ as illustrated in Fig. 3 and $Z^2(\bar{\mathbf{y}}) = \{(1, 3), (3, 8), (8, 14), (14, 24)\}$. We have illustrated a bipartite graph where each $j \in \mathcal{B}_\ell$ needs to be packed with items represented by a 3-tuple $(k, t_1, t_2)$ and each of these tuples have a weight of $t_2 - t_1$ equivalent to the number of hours in a given time block. The constraint (16b) implies that only one of the red arrow should be selected for $(k, t_1, t_2) = (1, 1, 3)$. The constraints (16c) ensure that for any given time period $t \in \mathcal{T}$, a battery $j$ is not assigned to more than one port. This is ensured in the formulation by summing over all the combinations of $(j, t_1, t_2)$ such that $t_1 + 1 \leq t \leq t_2$. This is also represented for $j = 2$ by green arrows and for $t = 4$. It implies that only one of the green arrows should be selected as $t = 4$ lies in the connected time blocks $(3, 8)$ and $(3, 6)$. Finally constraints (16d) and (16e) ensure that the demand for each battery $j \in \mathcal{B}_\ell$ is met as in the original model. The equivalence of (16) and original model(2) with fixing $\bar{\mathbf{y}}$ is apparent, and for brevity, we do not formally prove it.

**Observation 5.1** (Compactness of the Formulation). *The primary advantage of the bin-packing formulation (16) is its ability to significantly reduce the number of constraints and variables compared to the original model (2) with fixed $\mathbf{y} = \bar{\mathbf{y}}$. Specifically, the original model includes $B$ constraints for (2b)-(2c), $NT$ for (2d), $BT$ for (2e), $N(T-1)$ for (2f), and $BN(T-1)$ plus $T-1$ for (2g). Hence, the expression's growth is $O(BNT)$.*

*Analyzing the first constraint (16b) in model (16), we note that for any given $k \in \mathcal{N}$, there can be at most $T$ 2-tuples $(t_1, t_2) \in Z^k(\bar{\mathbf{y}})$, meaning the time horizon cannot be split into more than $T$ time-blocks. This estimate is conservative, as $\bar{\mathbf{y}} \in \tilde{P}^{\mathbf{y}\text{-s}}$ implies adherence to the switching constraint, potentially reducing the number of feasible time-blocks. Furthermore, the constraints for (16c) total $BT$, and for (16d)-(16e), a total of $B$ constraints exist, leading to an expression growth of $O(BT)$, given $B > N$. For the original model (2), which involves only the $\mathbf{x}$ variable, there are $BNT$ variables. In model(16) we have the same number of variable but assuming $T$ time-blocks for each $k \in \mathcal{N}$ in model (16) is overly conservative.*

In practice, for an instance with $B = 400, N = 140, \gamma = 45$, the number of variables was approximately 8,500 in model (16) compared to 1,000,000 for an equivalent model (2) with fixed $\mathbf{y}$. Similarly, the number of constraints was 189,880 versus 1,011,220. These numbers are before pre-solve, yet for (2), even pre-solving is time-consuming and negatively impacts computational efficiency. For larger instances, finding a feasible solution for the variable fixing problem using original model (2) often proves challenging, rendering it an ineffective Lagrangian heuristic unless our modification is considered.

### 5.2. $\sigma$-Local search using ergodic-sequences of Lagrangian subproblem solutions

In this section, we provide an overview of our approach to local search on an incumbent solution $\mathbf{x}^{\text{best}}$. For a predefined neighborhood parameter $\sigma$ where $N > \sigma > 0$, we identify a total of $\sigma$ ports. The batteries allocated to these ports will be relocated by moving them between the selected ports and time periods. The remaining parts of the solution, corresponding to the ports not included in the $\sigma$-neighborhood, are maintained as is. To select $\sigma$ ports we compute ergodic iterates, first presented in Larsson et al. (1999) which constructs approximations of the primal solution and the authors showed that the ergodic sequences in the limit converge to the optimal solution of the original problem. Later, Gustavsson et al. (2014) provided a more general technique to exploit subproblem solutions. Let us consider a sequence of ergodic iterates $\{\tilde{\mathbf{x}}^i\}$ constructed by using subproblem solutions $\mathbf{x}^*(\boldsymbol{\mu}^p), p < i$ obtained from solving the problem in $P^{\mathbf{x}}$ for multiplier $\boldsymbol{\mu}^p$ in previous iterations $p = 0, \ldots, i-1$ and defined as:

$$\tilde{\mathbf{x}}^i := \sum_{p=0}^{i-1} \phi_p^i \mathbf{x}^*(\boldsymbol{\mu}^p); \quad \sum_{p=0}^{i-1} \phi_p^i = 1; \quad \phi_p^i \geq 0, \quad p = 0, \ldots, i-1, \tag{17}$$

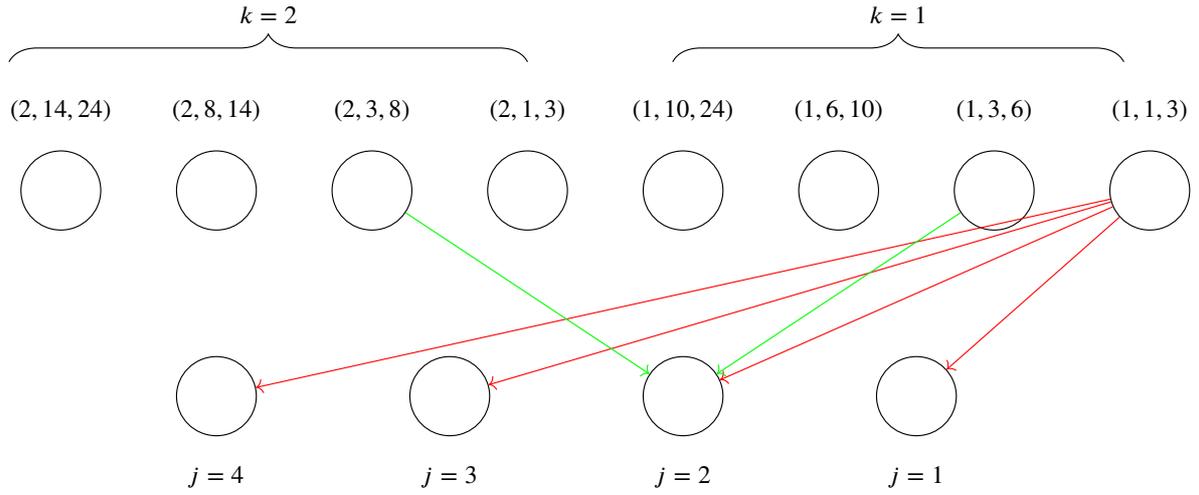

**Figure 7:** A bipartite graph representing the variable fixing problem. Each 3-tuple indicates $(k, t_1, t_2)$. The green arrows correspond to (16c) for $t = 3$, $j = 2$ and red arrows correspond to (16b) for $k = 1, (t_1, t_2) = (1, 3)$

where the convexity weights $\phi_p^i$ are chosen accroding to so-called $p^a$-rule:

$$\phi_p^i := \frac{(p+1)^a}{\sum_{r=0}^{i-1}(r+1)^a}, \quad p = 0, \ldots, i-1, i = 1, \ldots, ; \quad a \geq 0$$

Note that for higher values of $a$ the weights are shifted towards later subproblem solutions. It is important to note that since the vector $\mathbf{x}^*(\boldsymbol{\mu}^i)$ is a high dimensional vector storing it for all the iterations is impossible. To address that we utilize the concept of computing ergodic iterates as per (Gustavsson et al., 2014, eq. 17) (see (25) in Appendix 2 for more details ). Since we are taking convex combinations of subproblem solutions we can get non-binary values and hence, it remains to be seen how we utilize this information to obtain feasible solution.

We propose an approach that either results in a better solution than the incumbent or the same solution is returned. The core idea is to use ergodic iterate $\tilde{\mathbf{x}}^i$ which captures to a certain extent information from previous subproblem solutions (depending upon parameter $a$) to define the search space for a given parameter $\sigma$. Let us define an indicator function, $\chi_\sigma(k; \tilde{x}^i)$, which determines whether a given index $k$ belongs to the $\sigma$-neighborhood based on the total cost calculated from the ergodic iterate $\tilde{\mathbf{x}}^i$. Hence, each iteration we compute ergodic iterate and utilize them in a local search technique. We define the following indicator function:

$$\chi_\sigma(k; \tilde{x}^i) = \begin{cases} 1, & \text{if } \sum_{\ell \in \mathcal{L}} \sum_{j \in \mathcal{B}_\ell} \sum_{t \in \mathcal{T}} c_t \tilde{x}_{jkt}^i \text{ ranks among the top } \sigma \text{ highest sums for all } k' \in \mathcal{N}, \\ 0, & \text{otherwise.} \end{cases}$$

The $\sigma$-neighborhood, $\mathcal{N}_\sigma(\tilde{x}^i)$, is then defined as the set of all indices $k$ for which the indicator function returns 1, representing the indices with the highest aggregated contributions in terms of the ergodic iterate:

$$\mathcal{N}_\sigma(\tilde{x}^i) := \left\{ k \in \mathcal{N} \mid \chi_\sigma(k; \tilde{x}^i) = 1 \right\}; \quad \mathcal{B}_\ell(\mathbf{x}^{\text{best}}, \mathcal{N}_\sigma(\tilde{x}^i)) := \left\{ j \mid \exists t \in \mathcal{T}, x_{jkt}^{\text{best}} = 1, k \in \mathcal{N}_\sigma(\tilde{x}^i), j \in \mathcal{B}_\ell, \ell \in \mathcal{L} \right\}$$

The term $\mathcal{B}_\ell(\mathbf{x}^{\text{best}}, \mathcal{N}_\sigma(\tilde{x}^i))$ include all the batteries in the incumbent solution that are allocated to ports in $\mathcal{N}_\sigma(\tilde{x}^i)$. This allows us to selectively focus on those ports that have the highest cost across iterations. Note that if $\sigma = N$ then neighborhood includes all the ports, hence, the entire feasible set. Once we have identified $\mathcal{N}_\sigma(\tilde{x}^i)$ and $\mathcal{B}_\ell(\mathbf{x}^{\text{best}}, \mathcal{N}_\sigma(\tilde{x}^i))$ we can perform the local search on the incumbent $\mathbf{x}^{\text{best}}$ by solving the original model (2) with the exception that $\mathcal{B}_\ell$ needs to be replaced with $\mathcal{B}_\ell(\mathbf{x}^{\text{best}}, \mathcal{N}_\sigma(\tilde{x}^i))$ and similarly $\mathcal{N}$ with $\mathcal{N}_\sigma(\tilde{x}^i)$. Furthermore, certain constraints need to be adjusted to make sure that the resulting solution is still feasible. Hence, we either improve our best solution or we are returned the same solution. To be complete, we provide the model to be solved in the appendix (see (24)).

We provide a simple example in Fig. 8 to convey the main points. Let us consider $\mathbf{x}^{best}$ is known and illustrated in Fig. 8[left]. The rows are ports and columns are time periods. The text for each port and time period represent the allocated battery. Now for a given $\sigma$ and an ergodic iterate $\tilde{x}^i$, let us consider the set $\mathcal{N}_\sigma(\tilde{x}^i) = \{1, 2\}$. This implies only the first two ports are to be considered in the neighborhood. Hence, all the batteries and their respective hours can be reallocated to minimize the electricity cost as well as number of switching. Hence, batteries(hours) 1(1), 2(2), 4(3), 6(1), and 8(2) need to be reallocated between ports 1, 2 and also over all time periods. This is to be done by solving the optimization problem (24). The resulting solution with lower electricity cost and total number of switching is illustrated in Fig. 8[Right]. Note that the constraints in (24) ensure that the new solution is still feasible with the remaining partial solution that was not changed i.e. allocation done on ports 3 and 4. In the example we have used $\gamma = 4$. The reason why we use ergodic iterates and not the incumbent solution to define the neighborhood is by this way we can apply local search more frequently as ergodic iterate changes almost every iteration but same is not true for $\mathbf{x}^{best}$.

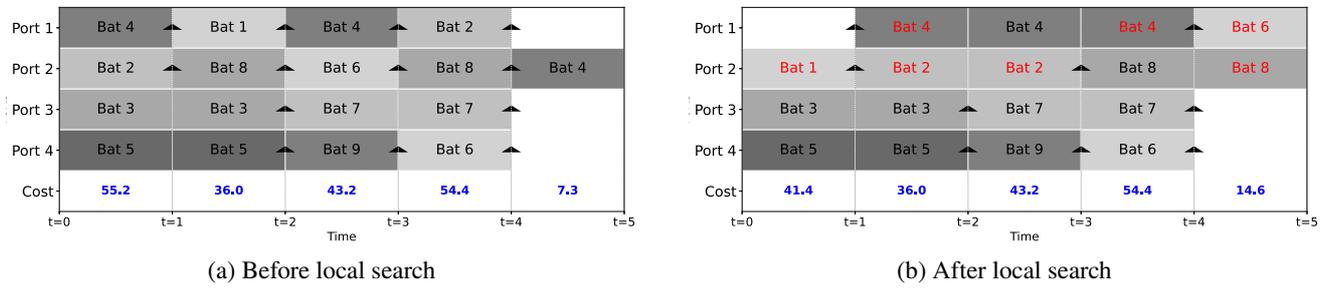

(a) Before local search

(b) After local search

**Figure 8:** 2-Local search on the first two ports. Left: Electricity cost: 196.1, Total switching: 13. Right: Electricity cost: 189.6, Total switching: 9. Cost (input) at the five time periods are [13.8, 9.0, 10.8, 13.6, 7.3] and $\gamma = 3$. The solid triangles represent a switching. The grey color coding is used simply to differentiate between batteries and the red text (Fig. 8b) are the batteries whose time period was changed w.r.t. the incumbent solution (Fig. 8a)

### 5.3. Warm starting the Lagrangian multipliers for a good initial lower bound

One of the principal challenges in subgradient optimization include initializing the Lagrangian multiplier, denoted as $\boldsymbol{\mu}^0$. Good initialization is important when dealing with a dual space containing a vast number of variables (i.e. dimension of $\boldsymbol{\mu}$), as is the case in our problem [e.g. for 200 batteries, $200 \cdot 100 \cdot 24 = 480,000$]. Due to the limitation on computation time, bad initialization can result in weaker lower bounds. Consequently, there is a need to develop so-called *warm start* strategy for the initial multiplier $\boldsymbol{\mu}^0$. We motivate our forthcoming suggestions by first demonstrating that any randomly chosen $\boldsymbol{\mu}^0$ results in weak lower bounds for our problem. For instance, consider initialization: $\mu_{jkt}^0 = 0$ or $\mu_{jkt}^0 = (1-\beta)c_t$. Notably, if $\mu_{jkt}^0 = 0$ for all $j, k, t$, it follows from equation (5) that $h_\mathbf{x}(\mathbf{0}) = \min_{\mathbf{x} \in P^\mathbf{x}} \sum_{\ell \in \mathcal{L}} \sum_{t \in \mathcal{T}_\ell} \sum_{j \in \mathcal{B}_\ell} \sum_{k \in \mathcal{N}} (\beta c_t) x_{jkt}$, and if $\mu_{jkt}^0 = (1-\beta)c_t$, the coefficient simplifies to just $c_t$ instead of $\beta c_t$. Both initialization lead to a low value of $h_\mathbf{x}(\boldsymbol{\mu}^0)$, as there are no penalties in $h_\mathbf{x}(\boldsymbol{\mu}^0)$ for the relaxed constraint (2f) to (2h). Hence, the optimal solution will exploit as many charging bays/ports as possible during time period $t$ with a low value of $c_t$. For the subproblem defined in $(\mathbf{y}, \mathbf{s})$ the objective of the Lagrangian subproblem is $\tilde{h}_{\mathbf{y},\mathbf{s}}(\boldsymbol{\mu}^0) = \min_{(\mathbf{y},\mathbf{s}) \in \tilde{P}^{\mathbf{y}\text{-s}}} \sum_{\ell \in \mathcal{L}} \sum_{t \in \mathcal{T}_\ell} \sum_{j \in \mathcal{B}_\ell} \sum_{k \in \mathcal{N}} ((1-\beta)c_t - \mu_{jkt})s_{jkt}$ (see (5)), a variable $s_{jkt}$ will be set to one only if there is a corresponding negative cost in the objective function i.e. $(1-\beta)c_t < \mu_{jkt}$, which is not the case for both initialization as $0 < \beta < 1$. In conclusion, if we use $\mu_{jkt}^0 = (1-\beta)c_t$, then $h_\mathbf{x}(\boldsymbol{\mu}^0)$ has equivalent optimal value as obtained by solving model(2) without constraints (2f)–(2h). This is a weak lower bound. Furthermore, $\tilde{h}_{\mathbf{y},\mathbf{s}}(\boldsymbol{\mu}^0) = 0$ because the objective coefficient of $s_{jkt}$ is zero (or positive when $\boldsymbol{\mu}^0 = 0$) and hence, none of the $y_{kt}$ will be set to one. Similar conclusion can be drawn for $\mu_{jkt}^0 = 0$.

This motivates us to explore better ways to initialize the Lagrangian multipliers. The general mechanism is as follows: at iteration $i$, the multiplier component $\mu_{jkt}^i$ increases if the subproblem solution component $x_{jkt}^*(\boldsymbol{\mu}^{i-1}) > s_{jkt}^*(\boldsymbol{\mu}^{i-1})$, decreases if $x_{jkt}^*(\boldsymbol{\mu}^{i-1}) < s_{jkt}^*(\boldsymbol{\mu}^{i-1})$, and remains unchanged if they are equal. Additionally, from the definition of $P^x$ (see (3a)), it is clear that $\mathbf{x}^*(\boldsymbol{\mu}^i)$ satisfies the demand for charging hours but fails to meet the switching constraint. In contrast, the subproblem solution $\mathbf{s}^*(\boldsymbol{\mu}^i)$ satisfies the switching constraint but not the demand. Since $h_{\mathbf{x}}(\boldsymbol{\mu}) = (\beta c_t + \mu_{jkt})x_{jkt}$, when $c_t$ is low, it is highly likely that $x^*_{jkt}(\boldsymbol{\mu}) = 1$. Consequently, periods with low costs are overused implying that for some time periods number of switching might be more that $\gamma$, leading to an increase in the corresponding multipliers. With this in mind, we applied subgradient optimization to a small instance with $|\mathcal{B}| = 20$, $|\mathcal{N}| = 10$, and $\gamma = 2$, noting that the relative dimensions of the three are similar to practical instances. Initially, for a few iterations the values of $\mu_{jkt}$ showed no apparent relation to the cost $c_t$ (see Figures 9a–9c), but after several iterations, as we approached the dual optimal value, a clear relationship emerged between the near-optimal Lagrangian multipliers and the cost (cf. Figures 9d–9f).

This is intuitive and simply indicates the fact that for the time periods when $c_t$ is low an optimal $\mu_{jkt}$ should be relatively higher as this should adjust the subproblem solution to the Lagrangian subproblem in $\mathbf{x}$ to encourage respecting the constraint (2h). Of course, due to demand windows and other input parameter there can be some deviations but for most of the time periods it should hold. We illustrate this by plotting a linear regression line for different iterations to understand the relation between $\mu_{jkt}$ and $c_t$ (Fig. 9). We indicate the final values of slope, intercept and coefficient of determinations as $\omega_1^\ell, \omega_2^\ell$ and $R^2, \ell \in \mathcal{L}$ respectively.

Our idea is to learn these parameters from small instances and then applying them to initialize $\boldsymbol{\mu}^0$ for larger, practical instances. This process does not require frequent recalculations. Instead, one can utilize previous day's prices to compute it well in advance. Empirical tests, conducted across a range of price fluctuations and instance variations, have demonstrated that the parameters typically remain stable without needing substantial modifications. Notably, near-optimal Lagrangian multipliers often exhibit a negative slope, accompanied by a significantly high coefficient of determination ($R^2$). Let us consider that these regression parameters are known then for any large-scale realistic instance an initial Lagrangian multiplier is defined as $\mu_{jkt}^0 = \hat{\mu}_{jkt} = \hat{\omega}_1^\ell c_t + \hat{\omega}_2^\ell, j \in \mathcal{B}_\ell, k \in \mathcal{N}, t \in \mathcal{T}_\ell, \ell \in \mathcal{L}$. In Tab. 3 we show the impact of using these warm start strategy by learning the regression parameter from a set of smaller instance i.e. $|\mathcal{B}| = 20$. We see the significant increase in the starting value of the lower bound.

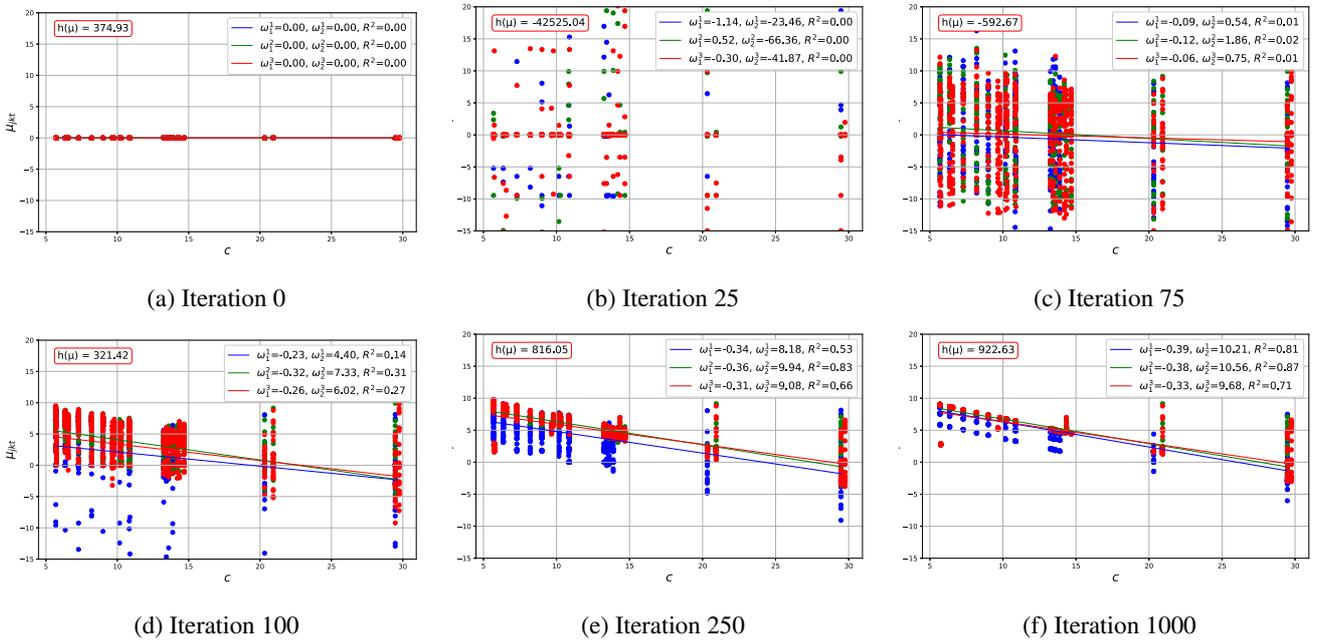

**Figure 9:** Evolution of the relation between $\boldsymbol{\mu}$ and cost $\mathbf{c}$ as we approach dual optimal solution for $B = 20, N = 10, \gamma = 2$. See a gif file: `https://bit.ly/3Yph9HD` [github]

## 5.4. Computing Lagrangian multipliers for enhanced lower bounds

The literature review highlights various methods for computing the subgradient, including a straightforward approach utilizing the formula discussed in the preceding section i.e. $d_{jkt}^i = g_{jkt}^i = (x^*(\boldsymbol{\mu}^i) - s^*(\boldsymbol{\mu}^i))$ for an iteration

**Table 3**
An overview of improved initial Lagrangian dual value due to warm starting $\mu^0 = \hat{\mu}$.

| $|\mathcal{B}|$ | $|\mathcal{N}|$ | $\gamma$ | $\tilde{h}(\mathbf{0})$ | $\tilde{h}(\hat{\mu})$ |
|---|---|---|---|---|
| 100 | 50 | 13 | 1843.80 | 4197.25 |
| 150 | 75 | 19 | 2768.26 | 6324.22 |
| 200 | 100 | 25 | 3713.29 | 8498.43 |
| 250 | 125 | 32 | 4663.45 | 10632.02 |
| 300 | 140 | 35 | 5687.09 | 13019.05 |
| 350 | 140 | 35 | 6859.21 | 15690.91 |
| 400 | 140 | 45 | 8148.44 | 17790.85 |

*i*. Alternatively, the modified deflected subgradient (MDS) method, as detailed by Belgacem and Amir (2018) (refer to the Literature section or eq. (1)), can be employed. While both methods exhibit convergence properties, executing numerous iterations becomes impractical for large-scale instances.

Moreover, we encounter a distinctive phenomenon: irrespective of the subgradient update method employed or step length parameter, the plot for the lower bounds invariably exhibits a sharp decline followed by a gradual increase, as illustrated in Fig. 10[7]. This leads to numerous iterations being expended without yielding any enhancement in the best-known lower bounds. This phenomenon arises due to the fact that, in the initial iterations, the Lagrangian multipliers $\mu$ are not sufficiently large, leading to numerous cases where $x_{jkt} > s_{jkt}$. However, after several iterations of consistent increase in the Lagrangian multipliers, some $s_{jkt}$ values begin to reach one, and it may occur that $s_{jkt} > x_{jkt}$. This scenario can result in $\tilde{h}_{\mathbf{y},\mathbf{s}}(\mu^i)$ becoming negative when $\mu^i_{jkt} \gg (1-\beta)c_t$ (refer to Eq. (5)). It requires a significant number of iterations to obtain multipliers that achieve $x_{jkt} = s_{jkt}$. This phenomenon exists even if we do not use the warm starting strategy and simply initialize with zero. In order to mitigate this we propose an intermediate step. First, we calculate an intermediate Lagrangian multiplier after the $i^{th}$ iteration as follows:

$$\mu_{jkt}^{i+\frac{1}{2}} := \mu_{jkt}^i + \tau g_{jkt}^i, \qquad \forall j,k,t, \tag{18a}$$

$$\mu_{jkt}^{i+1} := \omega_1^\ell(\boldsymbol{\mu}^{i+\frac{1}{2}})c_t + \omega_2^\ell(\boldsymbol{\mu}^{i+\frac{1}{2}}), \qquad \forall j \in \mathcal{B}_\ell, k, t, \ell, \tag{18b}$$

where $\omega_1^\ell(\boldsymbol{\mu}^{i+\frac{1}{2}})$ and $\omega_2^\ell(\boldsymbol{\mu}^{i+\frac{1}{2}})$ represent the linear regression slope and intercept, respectively, calculated using $\boldsymbol{\mu}^{i+\frac{1}{2}}$ and $\mathbf{c}$. The subgradient, can be derived either through MDS or simply by calculating the difference between $\mathbf{x}$ and $\mathbf{s}$, as previously discussed. The resulting Lagrangian dual, $\tilde{h}(\boldsymbol{\mu}^{i+1})$, remains a valid lower bound, as there are no constraints on $\boldsymbol{\mu}$. However, this approach forfeits explicit guarantees on convergence rates, leaving it as an open question whether reasonable convergence bounds can be achieved. Hence, it should be considered as a limitation. Nonetheless, the forthcoming computational results appear very promising in improving lower bounds and hence, lower optimality gap.

To demonstrate the variations in the Lagrangian dual and regression parameters both with and without the implementation of the described update rule, we introduce a comparative analysis depicted in Fig. 11 for the instance $\mathcal{B} = 200$. On the right side of the figure, a notable decrease in the Lagrangian dual (lower bound) is observed when projections onto the regression line are not used to update Lagrangian multipliers. In contrast, the left side of the figure showcases a modest increase in the Lagrangian dual when described projections are considered. We differentiate between two approaches: *Reg*, indicating that projection onto the regression line is performed, and *nReg*, where no such projection is made. Following this classification, we further describe the method employed for computing the subgradient. We differentiate between *Simple*, where the subgradient $d_{jkt}^i = g_{jkt}^i$ is computed as the difference between $x_{jkt}^*(\mu^i)$ and $s_{jkt}^*(\mu^i)$, and *MDS*, which utilizes the modified deflected subgradient method as detailed in (1). Note that even though it seems clear right now that these projections help in improving lower bounds but it may have positive or negative impact on upper bounds which is yet to be tested. We do a detailed empirical testing in the results section.

---

[7]It should be noted that a fixed time limit is applied across all instances. Consequently, as the size of an instance increases, the number of iteration reduces.

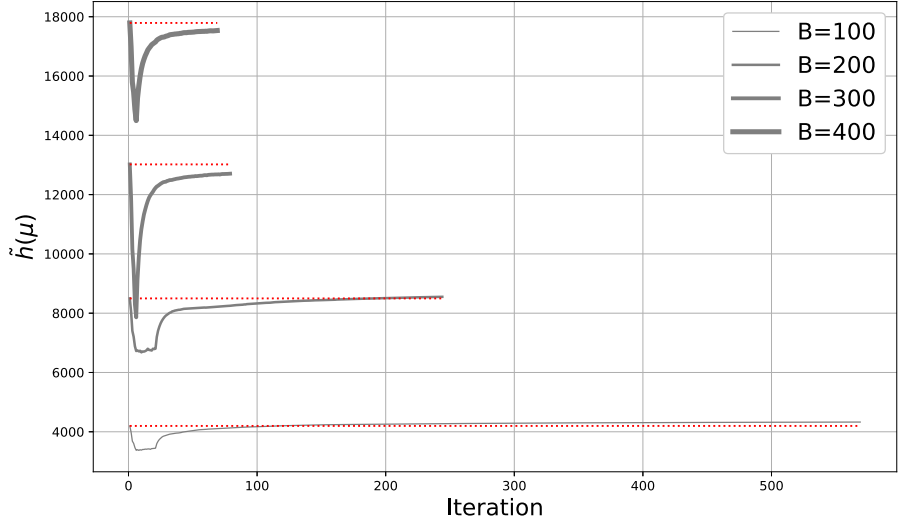

Figure 10: Reduction in Lagrangian dual value for different instances using standard Subgradient update method. The red dotted line is drawn at $\tilde{h}(\hat{\boldsymbol{\mu}})$

### 5.5. Summary of the algorithm

A summary of the overall algorithm is provided in Alg. 1. We start with initializing step length parameters, lower and upper bounds. Furthermore, we initialize the multiplier as $\boldsymbol{\mu}^0$ as indicated in Sec. 5.3. We initiate a while loop that will be terminated if time limit exceeds or we reach $\epsilon$-optimality gap. The first step is to solve the Lagrangian subproblems corresponding to $\mathbf{x}$ and $\mathbf{y}, \mathbf{s}$. The former to be solved using Google OR's min-cost network flow library (see Prop. 4.3) and the latter is be solved using the time-block formulations as indicated in Sec. 4.1. Both these subproblems are to be solved in parallel. Furthermore, if we have $Loc$ = True i.e. we are using local search then we check the condition on line 5 which is True after every $\frac{heu}{2}$ iteration. Note that it is neither realistic to apply variable fixing heuristic nor local search in every iteration as solving the model takes some time and generally the difference in solutions is not significant in each iteration. Hence, both variable fixing and heuristics are applied after a certain number of iterations have passed[8]. Note that the local search can also be solved in parallel to the two Lagrangian subproblems. Hence, from lines 4-7 the elapsed time is equal to the problem which takes most time. Each of these problems have a time limit of 50 seconds. On lines 9-14 depending upon the type of subgradient update rule we compute the direction $\mathbf{d}^i$. On lines 15-19 we apply variable fixing heuristic using bin-packing formulation as indicated in (16). It is followed by updating the best solutions and step length parameter to update the Lagrangian multiplier for the next iteration. The update is done depending upon whether we project on the computed regression line as discussed in (18a)-(18b).

## 6. Computational experiments

We provide an overview of the input data, the computational setup, and the experimental design. All the algorithms are written using Spyder IDE for Windows system(s) with a 1.70 GHz processor, 16 GB RAM, and 4 cores. We use `Python 3.9` as the interpreter and `Gurobi 9` as the solver. The sub-problems in $\mathbf{x}$ and $(\mathbf{y}, \mathbf{s})$ and local search model are solved on parallel identical computers in a cluster. Each computer solves a sub-problem and its full processor capacity and all CPU cores are available when solving a given sub-problem.

### 6.1. Instances

We assume that we have homogeneous chargers with an output of 100W. Given that the maximum output power at a facility can reach 14,000W, no more than 140 ports can be utilized simultaneously. We have based our calculations on a typical intra-day electricity price profile expressed in öre per 100Wh (where 100 öre equals 1 SEK). Although all batteries are homogeneous, they can have varying states of charge (SoC). Based on the SoC, we assume that users of our model apply a rule to compute input sets $\mathcal{B}_\ell, \ell \in 1, 2, 3$. This may depend on the demand prediction algorithms used by the operator to adjust the number of batteries required at different times of the day.

---

[8] we use the parameter heu for variable fixing and $\lfloor \frac{heu}{2} \rfloor$ for local search

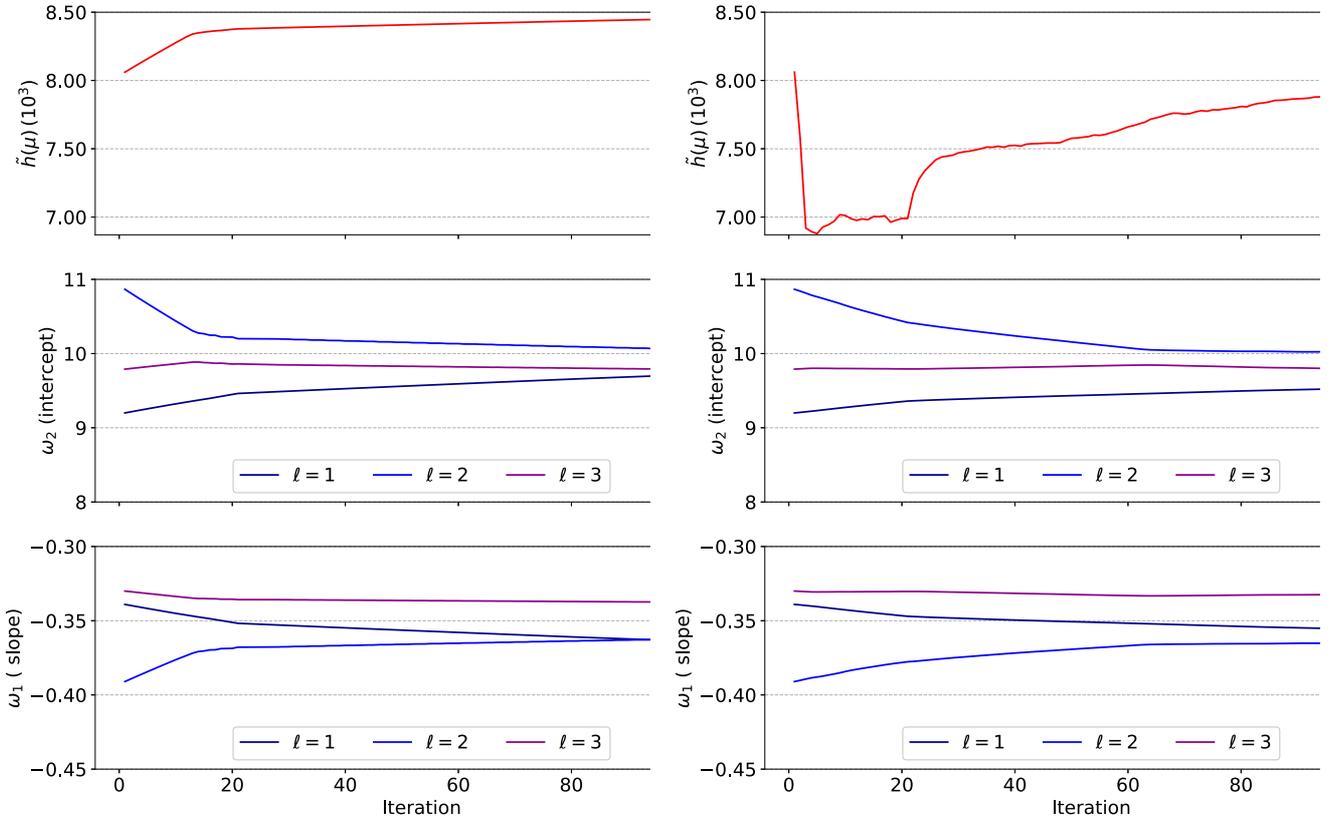

**Figure 11:** Stacked plot to illustrate impact of different update strategy for the instance with 200 batteries. Left(Reg): projection onto the regression line and Right(nReg): No projections.

For staff comfort, switchings are permitted only after a minimum interval of one hour. This interval can be extended or even adjusted to a variable discretization of the planning horizon, depending on specific requirements. The model is designed to be sufficiently flexible to accommodate various modifications without compromising performance stability. The test instances are available here: https://bit.ly/3EHxssF [github]. Initially, we test our algorithms with electricity price as illustrated in Fig. 18a with standard peak price at time 08:00 and 20:00. We refer to this as the *base* case.

### 6.2. Computational benchmarking

This section is organized to clearly distinguish between our proposed algorithm and other state-of-the-art approaches. Additionally, a separate subsection explores the sensitivity of results from our algorithm when varying certain components. This analysis provides detailed insights into individual contributions towards improvements in both the lower and upper bounds. This is important since, in exact methods, introducing components to enhance the upper bounds may negatively affect the lower bounds, and vice versa. The key is to find a good balance based on user requirements.

**Table 4**
Comparison of our algorithm SOVLR with other solvers under computation time limit $t^1 = 1800$ secs.

| $|\mathcal{B}|$ | $|\mathcal{N}|$ | $\gamma$ | $\bar{h}^{\text{best}}$ | | | $\underline{h}^{\text{best}}$ | | | Gap[%] | | |
|---|---|---|---|---|---|---|---|---|---|---|---|
| | | | SOVLR | GRB | CP | SOVLR | GRB | CP | SOVLR | GRB | CP |
| 100 | 50 | 13 | 4494.10 | 4975.40 | 4537.81 | 4264.13 | 4088.52 | 3687.67 | 5.11 | 17.80 | 18.33 |
| 150 | 75 | 19 | 6749.24 | 8613.10 | 7125.63 | 6418.94 | 5536.70 | 5207.29 | 4.89 | 35.71 | 26.92 |
| 200 | 100 | 25 | 9077.46 | 11482.22 | 9680.50 | 8589.59 | 7426.70 | 6698.53 | 5.37 | 35.90 | 30.80 |
| 250 | 125 | 32 | 11342.48 | 14246.05 | 12317.56 | 10725.49 | 9327.10 | 7611.05 | 5.43 | 34.52 | 38.20 |
| 300 | 140 | 35 | 13790.88 | - | 15960.09 | 13115.81 | 11374.26 | 8658.11 | 4.89 | 100 | 45.75 |
| 350 | 140 | 35 | 16723.77 | - | - | 15786.15 | 13718.49 | - | 5.60 | 100 | - |
| 400 | 140 | 45 | 18760.60 | - | - | 17878.65 | 16296.94 | - | 4.70 | 100 | - |

**Table 5**
Comparison of algorithms under computation time limit $t^2 = 3600$ secs.

| $|\mathcal{B}|$ | $|\mathcal{N}|$ | $\gamma$ | $\bar{h}^{\text{best}}$ | | | $\underline{h}^{\text{best}}$ | | | Gap[%] | | |
|---|---|---|---|---|---|---|---|---|---|---|---|
| | | | SOVLR | GRB | CP | SOVLR | GRB | CP | SOVLR | GRB | CP |
| 100 | 50 | 13 | 4475.76 | 4700.32 | 4515.26 | 4274.60 | 4259.74 | 3687.67 | 4.40 | 9.37 | 18.32 |
| 150 | 75 | 19 | 6749.24 | 8613.10 | 7000.25 | 6425.30 | 5536.70 | 5236.23 | 4.79 | 35.71 | 25.19 |
| 200 | 100 | 25 | 9057.61 | 11482.22 | 9490.60 | 8611.47 | 7426.70 | 6712.08 | 4.92 | 35.90 | 29.27 |
| 250 | 125 | 32 | 11342.48 | 14246.05 | 11950.39 | 10737.39 | 9327.10 | 7743.16 | 5.33 | 34.52 | 35.20 |
| 300 | 140 | 35 | 13790.88 | 17716.19 | 14826.85 | 13131.72 | 11374.26 | 8715.56 | 4.77 | 35.79 | 41.22 |
| 350 | 140 | 35 | 16723.76 | - | - | 15798.45 | 13718.49 | - | 5.53 | 100 | - |
| 400 | 140 | 45 | 18760.60 | - | - | 17898.33 | 16296.94 | - | 4.59 | 100 | - |

### 6.2.1. Comparison with other state-of-the art methods

We have explored three distinct solution procedures, including our own, referred to as SOVLR (Alg. 1), Gurobi's MIP solver (GRB), and Google's constraint programming SAT solver (CP). The choice of GRB and CP was driven by the fact that the former is one of the best commercially available MIP solvers, whereas the latter is believed to be more adept at handling the logical constraints present in our model. We tested two different time limitations to observe the benefits of extending the allowed solution time. Initially, in Table 4, we imposed a time limit of 1800 seconds and noted that both GRB and CP failed to compute feasible solutions for instances with 350 and 400 batteries. It is noteworthy that the SOVLR algorithm incorporates all suggested improvements from the previous section, except for local search. We will further discuss the variations in the algorithmic components and clarify their impact on the lower and upper bounds. Comparing results in Tables 4 and 5 reveals that increasing the solution time by 1800 seconds slightly enhanced the results, though not significantly. Therefore, a user with limited time can still achieve good results with an 1800-second limit. Additionally, only our algorithm consistently managed to compute reasonable feasible solutions across all instances. Interestingly, even though the instance with 400 batteries contains more variables, it exhibits a lower optimality gap than the instance with 350 batteries. This is primarily because we increased the value of $\gamma$ from 35 to 45, thereby relaxing the problem. It was concluded that an instance with 400 batteries and $\gamma = 35$ is not feasible. In all instances, our proposed algorithm managed to achieve an optimality gap of less than 6%, which is satisfactory considering the application and time constraints. It is also interesting to note that GRB manages to compute better lower bound than CP but for upper bounds reverse is true i.e. CP was much better than GRB. Now, in the upcoming section we focus our attention to understanding SOVLR in more detail and also trying some perturbation on the electricity price to see how it impacts the relative performance.

### 6.2.2. Comparison of different variations of SOVLR

In the preceding sections, we have presented several techniques to improve the solving algorithm. This subsection presents results from benchmarking tests to uncover the impact of these variations on lower bounds, upper bounds, the optimality gap, and primal-dual integrals (definitions to follow). In Sec. 4.1, we presented the time-block reformulation for the problem in $(\mathbf{y}, \mathbf{s})$. In Fig. 12, we illustrate the benefits of using the time-block on the upper and lower bounds. Clearly the reformulation has higher lower bound and lower upper bound as compared to the original formulation.

Next, we present the impact of variations on our algorithm. In Table 6, we provide an overview of different options available for algorithmic components. Firstly, we discuss the method, which indicates how we compute subgradients. Two options are available, as indicated in Sec. 5.4 §1: *Simple* and *MDS*. After computing the subgradient, we need to update the Lagrangian multiplier. For this purpose, we used an intermediate step described in Sec. 5.4 §2. We also include the option of omitting the intermediate regression step, indicated as *nReg*. Next, we present two options for handling the partial-variable fixing problem: using the bin packing formulation (*var-fix-binP*), as described in Sec. 5.1, or simply fixing the variables (*var-fix*). Finally, we have the option to use local search, indicated by *Loc*, as described in Sec. 5.2. One might question why we consider variations without these contributions if they yield good results. The answer is that each contribution, while beneficial, may consume time and reduce the number of iterations, potentially having an overall negative impact. Hence, it is important to understand the impact of each variation after we integrate all the components.

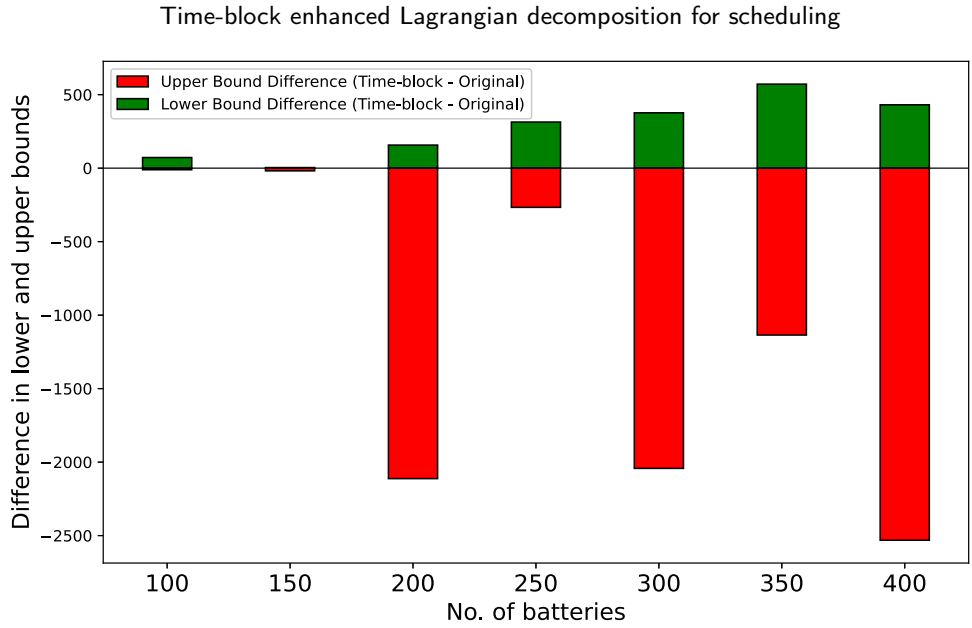

**Figure 12:** Reduction in upper bounds and increase in lower bounds by using time-block reformulation for the Lagrangian sub-problem in $\mathbf{y}, \mathbf{s}$ as described in Sec. 4.1 as opposed to the original.

**Table 6**
Attributes included in different variations of our proposed algorithm are marked in black.

| | | var. | 1 | 2 | 3 | 4 | 5 | 6 | 7 | 8 | 9 | 10 | 11 | 12 | 13 | 14 | 15 | 16 |
|---|---|---|---|---|---|---|---|---|---|---|---|---|---|---|---|---|---|---|
| Method | Sec. 5.4 §1 | Simple | | | | | | | | | | | | | | | | |
| | | MDS | | | | | | | | | | | | | | | | |
| Update | Sec. 5.4 §2 | Reg | | | | | | | | | | | | | | | | |
| | | nReg | | | | | | | | | | | | | | | | |
| Heuristic | Sec. 5.1 | var-fix-binP | | | | | | | | | | | | | | | | |
| | | var-fix | | | | | | | | | | | | | | | | |
| Local | Sec. 5.2 | Loc | | | | | | | | | | | | | | | | |
| | | nLocal | | | | | | | | | | | | | | | | |

We have encoded these variation with numbers in the range $\{1, \ldots, 16\}$ as indicated in column 4 onwards in Tab. 6. Next, we need good performance measures to benchmark these results. We must also think that different users may have different time limitations, hence, a performance indicator should not be too reliant on our time limitations. Hence, we utilize lower and upper bounds as well as both the optimality gap and the primal-dual integral. The latter is useful for understanding the speed of convergence, as two variations that achieve the same lower and upper bounds may differ in how quickly they approach the obtained solution. For instance, in Figure 13b, we illustrate the convergence plot of two variations, namely 5 and 6, of our algorithm and observe that they obtain similar lower and upper bounds. However, in this case, the primal-dual integral of variation 5 is lower—indicating a better performance. This is also helpful for the end-users having less available computational time and then they can focus on primal-dual integral. In Figure 13a, which compares variations 2 and 6 for an instance involving 350 batteries, the convergence plot appears quite different, yet the primal-dual integral is same. We formally define the primal-dual integral as follows: The primal-dual integral

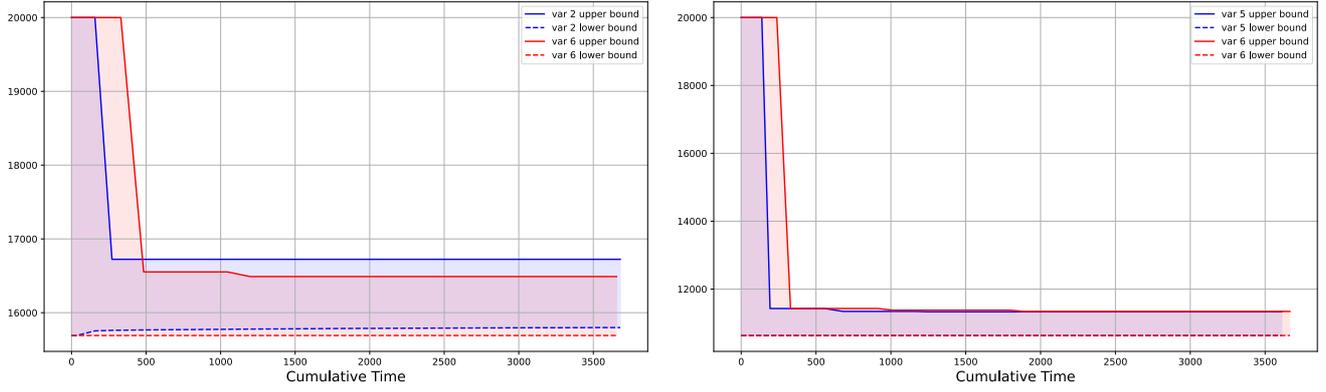

(a) Primal-dual integral for instance with 350 batteries and comparison with variations 2 and 6

(b) Primal-dual integral for instance with 250 batteries and comparison with variations 5 and 6

**Figure 13:** Comparison of primal-dual integrals

$\Theta$ over a given time limit $t^{\text{lim}}$ is defined as:

$$\Theta(t^{\text{lim}}) = \sum_{i=1}^{t^{\text{lim}}} \left(\bar{h}(t_i) - \underline{h}(t_i)\right) \cdot (t_i - t_{i-1}),$$

where $\bar{h}(t_i)$ represents the best known upper bound on the optimal solution value at discrete time $t_i$, and $\underline{h}(t_i)$ denotes the best known lower bound at the same point. The term $t_i$ indicates discrete time points, with $t_0$ as the initial time (often zero). The term $t_i - t_{i-1}$ measures the duration of the interval between consecutive time points $i$ and $i-1$. Thus, $\Theta$ quantifies the area between the upper and lower bounds over the time range from $t_0$ to $t^{\text{lim}}$.

First, we present a heatmap in Figure 15, summarizing the upper bounds achieved by different variations. The term "NA" indicates that variations 3, 4, 7, 8, 11, 12, 15, and 16 did not find any feasible solution in larger instances. These variations exclude our suggested bin-packing formulation, thereby highlighting the benefits of our heuristic based on this approach. We have removed these variations from further analysis to improve visualization quality and avoid skewing results with less effective variations.

Fig. 16 illustrate heatmaps for lower and upper bounds. The color coding is adjusted depending upon the normalization used for each type of instance. For example, for the instance with 100 batteries the variation 14 resulted in highest lower bound. Hence, its color code correspond to the value 0, where as for other variations we compute the percentage decrease as compared to the best-known lower bound. Similarly, for the upper bound, we calculate the percentage increase compared to the variation with the least upper bound for each iteration. Mean and median (in bracket) values are mentioned at the top of each variation column. These summary statistics are for the percentage increase of upper bounds and decrease for lower bounds. According to the median values the variation 10 (which uses MDS, Reg and no Local search) has the lowest value of 0.08 for lower bound. It seems to be logical as using local search generally improves the upper bound but since it takes away some time from the time allowance it may result in reducing the number of iterations. Consequently, corresponding lower bound can be less and optimality gap higher. In Fig. 16b, the best performance as per median values is var= 13 for upper bounds. This variation uses MDS, nReg and local search. Again, it is expected that the variation with local search will have slight advantage in computing better upper bounds. Mean and median values (the latter in brackets) are displayed at the top of each column for each variation. In Figure 16b, the best performance according to median values is from variation 13 for upper bounds. This variation utilizes MDS, nReg, and local search. This result seems logical since employing local search typically improves the upper bound. However, because it consumes some of the available time, it may reduce the number of iterations. Consequently, the corresponding lower bound may be lower and the optimality gap higher.

In Figure 17a, we illustrate the optimality gap (expressed as a percentage). Since the optimality gap is already normalized, we have provided the actual mean and median values for each variation. Note that the optimality gap is an effective measure to assess the balance between the upper and lower bounds. Variations 2, 9, and 10 all display reasonable values. Specifically, the second variation, which employs Simple, Reg, and no local search, consistently

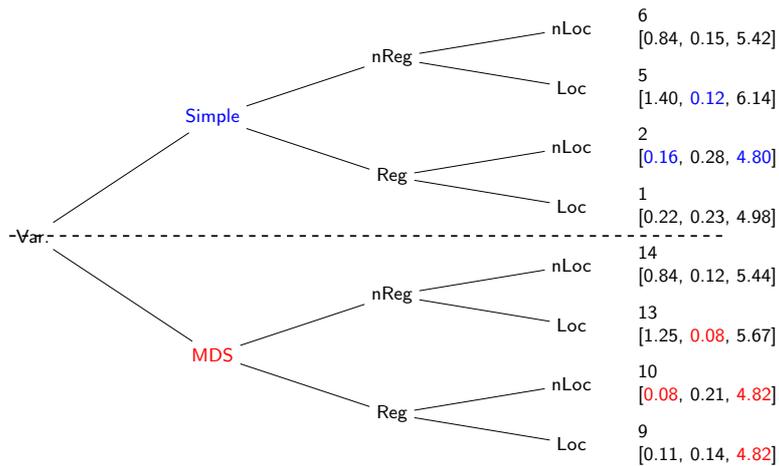

**Figure 14:** Details of different algorithmic variations. The median values are of [lower bound (normalized), upper bound (normalized), optimality gap (percentage)]. All the variations are using time-block reformulation (model 12) and our bin-packing heuristic (16)

shows an optimality gap of less than 6% across all instances. In Figure 17b, we demonstrate the primal-dual integral (pdi), normalized by dividing by the best pdi for a given instance. According to the column medians, variations 9 and 10 appear to be the most effective.

It is also valuable to summarize these medians in a tree structure as shown in Figure 14. The initial decision concerns the choice of subgradient method. For Simple, the best median values for both the lower bound and optimality gap are achieved using Reg without local search (variation 2). Conversely, as anticipated, the best upper bounds are obtained when local search is combined with nReg (variation 5). From our analysis, it appears that the subproblems are slightly more challenging to solve when multipliers are projected onto the regression line. Similarly, when MDS is employed, the best results for upper bounds are achieved using nReg and local search (variation 13), while the best values for lower bounds and optimality gap occur with Reg and no local search (variation 10). Overall, if the focus is solely on upper bounds, variation 13 may be beneficial, as it also shows strong performance in terms of the upper bound and the primal-dual integral, as illustrated in Figure 17b. However, both variations 2 and 10 are well-balanced in terms of both lower bounds and optimality gaps. Finally, in Fig. 22 we illustrate a gantt chart obtained from the solution of our algorithm for the second variation for the instance with 100 batteries.

Hence, variations 2, 5, 10, and 13 are interesting options. Furthermore, we also test these variations on a slightly perturbed electricity price profile (see 18b). It is common during summer for peak hours to be extended, as indicated in Fig. 18b. We now summarize the results for different price profiles, both base and extended. In Table 7, we present results for the variations of interest across two different price cases. In columns 5-8, we utilize different variations and report their respective upper bounds. In the last column, we use the best-known lower bound to calculate the actual optimality gap in columns 10-14. In column 9, we report the best upper bound from either GRB or CP. Interestingly, the optimality gaps generally increase with the new extended case. Furthermore, neither CP nor GRB could identify feasible solutions for instances with 250 or more batteries. Surprisingly, variations 2 and 10 did not identify a feasible solution for the instance with 350 batteries, and variation 2 also failed to find a feasible solution for the instance with 400 batteries. This is consistent with Fig. 14; since we do not use local search and using *Reg* generally results in worse upper bounds. Since our user is interested in upper bounds, using variations 5 and 13 is prudent. Both seem to perform similarly, as per Table 7. Both these variations use our time-block reformulation for the subproblem in $(\mathbf{y}, \mathbf{s})$ (12), Lagrangian heuristic using bin-packing formulation (16), and local search approach using ergodic-iterates as described in Sec. 5.2. Furthermore, it is extremely useful to use the initialized lower bound as per Sec. 5.3. The only suggested contribution that is not considered in variations 5 and 13 is from Sec. 5.4. Although variations which utilize proposed intermediate projection step provide better lower bounds but the upper bounds i.e. feasible objective values are not good Fig. 16 and see variations 2 and 10 in Tab. 7 both include Reg (see Tab. 6).

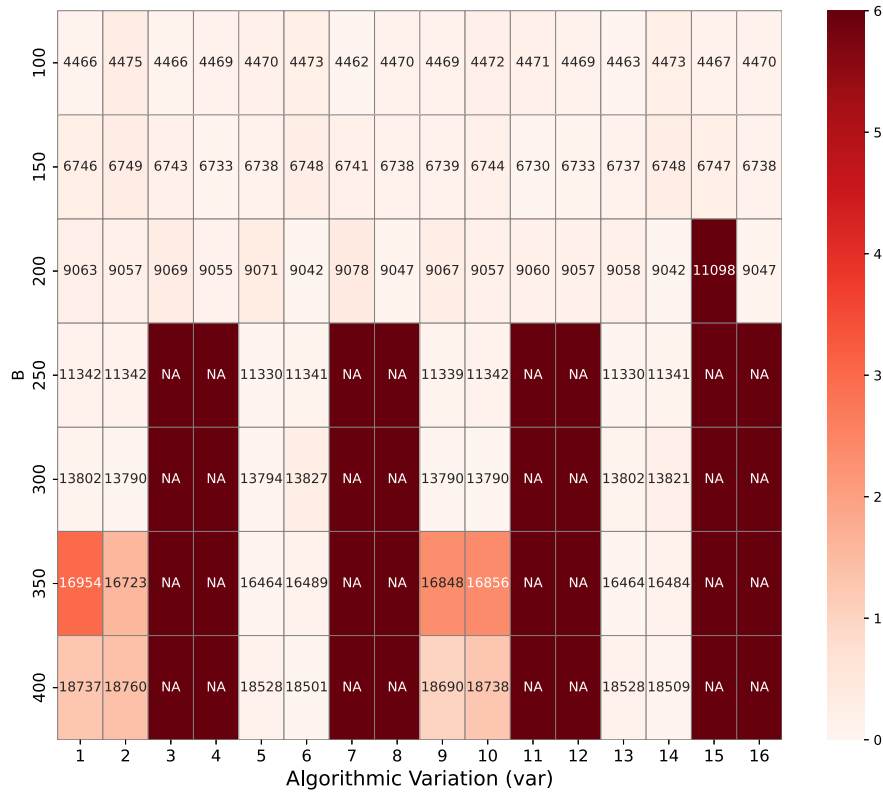

**Figure 15:** Impact on upper bounds of algorithmic variations. NA imlies no feasible solution was found which happens for variations with *var-fix* heuristic i.e. when bin-packing formulation is not used

**Table 7**
Final results table with two different peak price case i.e. base and extended. Time limit used is 3600 seconds.

| B | N | $\gamma$ | case | $\bar{h}$ (variations) | | | | | Gap% | | | | | best $\underline{h}$ |
|---|---|---|---|---|---|---|---|---|---|---|---|---|---|---|
| | | | | 2 | 5 | 10 | 13 | CP | 2 | 5 | 10 | 13 | CP | |
| 100 | 50 | 13 | base | 4475.8 | 4470.3 | 4472.8 | 4463.9 | 4515.3 | 2.9 | 2.7 | 2.8 | 2.7 | 3.7 | 4345.7 |
| 150 | 75 | 19 | base | 6749.3 | 6738.2 | 6744.0 | 6738.0 | 7000.2 | 4.7 | 4.5 | 4.6 | 4.5 | 8.1 | 6433.3 |
| 200 | 100 | 25 | base | 9057.6 | 9071.8 | 9057.6 | 9059.0 | 9490.6 | 4.8 | 5.0 | 4.8 | 4.9 | 9.2 | 8619.2 |
| 250 | 125 | 32 | base | 11342.5 | 11330.7 | 11342.5 | 11330.7 | 11950.4 | 4.8 | 4.7 | 4.8 | 4.7 | 9.6 | 10797.3 |
| 300 | 140 | 35 | base | 13790.9 | 13794.6 | 13790.9 | 13802.2 | 14826.8 | 4.4 | 4.4 | 4.4 | 4.5 | 11.1 | 13183.5 |
| 350 | 140 | 35 | base | 16723.8 | 16464.1 | 16856.8 | 16464.1 | - | 5.4 | 3.9 | 6.1 | 3.9 | - | 15823.8 |
| 400 | 140 | 45 | base | 18760.6 | 18528.6 | 18738.0 | 18528.6 | - | 4.6 | 3.4 | 4.5 | 3.4 | - | 17898.3 |
| 100 | 50 | 13 | ext. | 4697.2 | 4649.0 | 4667.7 | 4649.0 | 4904.1 | 5.7 | 4.8 | 5.1 | 4.8 | 9.7 | 4427.7 |
| 150 | 75 | 19 | ext. | 7234.6 | 7001.9 | 7200.3 | 7001.9 | 9404.2 | 8.3 | 5.3 | 7.9 | 5.3 | 29.4 | 6633.7 |
| 200 | 100 | 25 | ext. | 9739.6 | 9423.7 | 9732.4 | 9422.2 | 13670.7 | 9.1 | 6.0 | 9.0 | 6.0 | 35.2 | 8855.5 |
| 250 | 125 | 32 | ext. | 12124.5 | 11770.7 | 12057.7 | 11770.7 | 15897.6 | 9.4 | 6.6 | 8.9 | 6.6 | 30.8 | 10990.8 |
| 300 | 140 | 35 | ext. | 14973.4 | 14310.0 | 14956.8 | 14313.0 | - | 10.8 | 6.7 | 10.7 | 6.7 | - | 13354.8 |
| 350 | 140 | 35 | ext. | - | 18169.0 | - | 18249.7 | - | - | 12.0 | - | 12.4 | - | 15981.8 |
| 400 | 140 | 45 | ext. | - | 20212.8 | 24645.5 | 20216.2 | - | - | 9.0 | 25.3 | 9.0 | - | 18400.0 |

## 7. Sensitivity analysis

Now we analyze the impact of varying parameters such as $\gamma$ (limitation on the number of switchings) and $N$ (# of ports) on cost with fixed number of batteries. Furthermore, we also analyze how different electricity price plans may impact the cost and number of switchings for an operator.

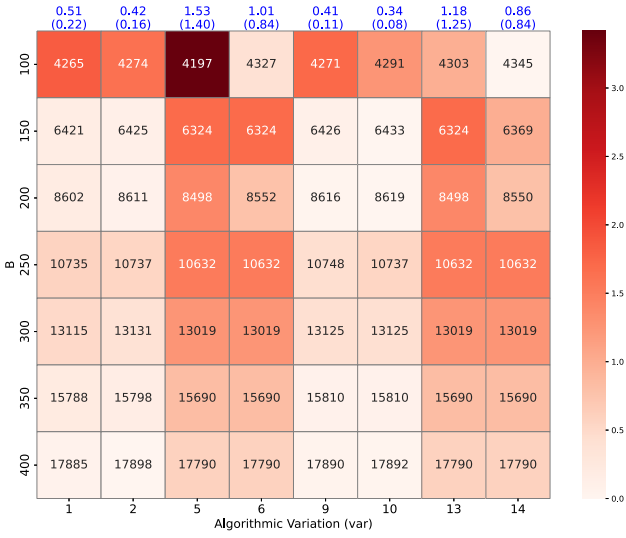

(a) Best lower bounds obtained

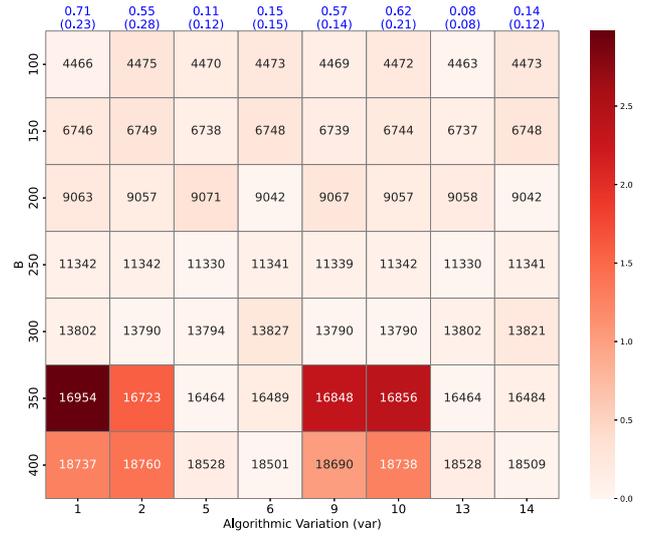

(b) Best upper bounds obtained

**Figure 16:** Comparison of different variations w.r.t. lower and upper bounds

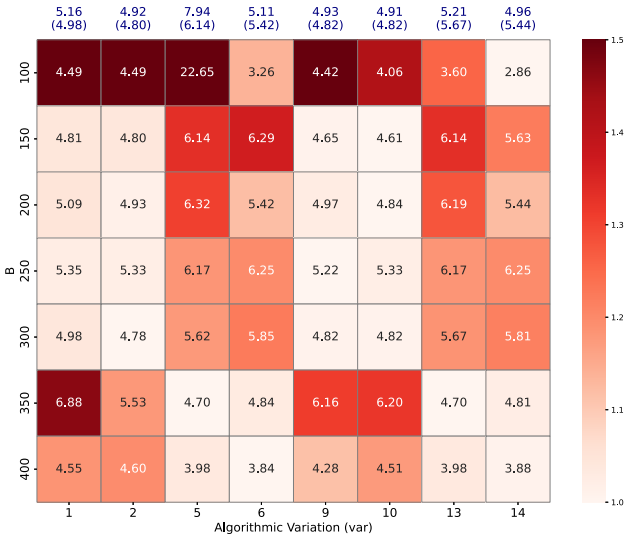

(a) optimality gap

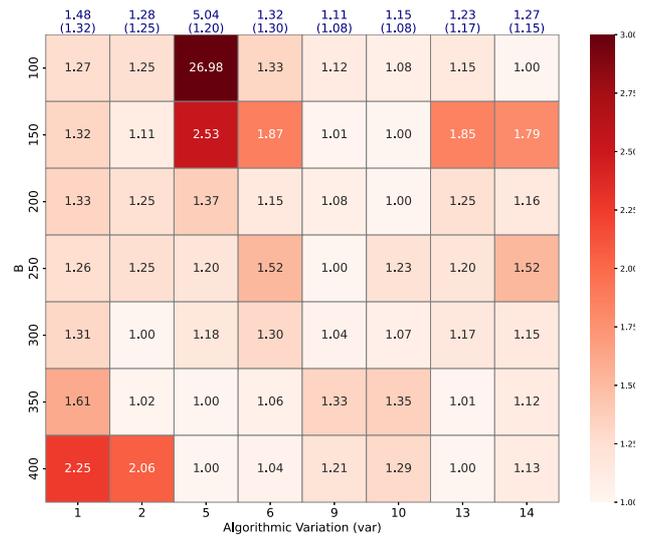

(b) primal-dual integral

**Figure 17:** Comparison of different variations w.r.t. lower and upper bounds

Firstly, in Figures 19a and 19b, we summarize the results obtained by varying $\gamma$ and the number of ports $N$ across two instances: one with 100 batteries and the other with 300 batteries, respectively. In Fig. 19a, it is obvious that higher value of $\gamma$ increases (relaxes) the feasible set as compared to lower values of $\gamma$, and hence, we obtain a lower objective value. However, this difference/benefit of having a higher value of $\gamma$ diminishes as the number of ports reduces. In Fig. 19b, we repeat the same experiment but for the instance with $B = 300$. Both plots are generated using variation 13 as our focus was on getting the best feasible solution. It is a good idea for the user to develop similar plots to assess the trade-off between the electricity cost and number of ports ($N$) for their instance to analyze best possible option for $\gamma$ and $N$.

It is also common for several commercial charging facilities to purchase different pricing plans, which may impact the number of switchings and the total electricity cost. In this section, we introduce three possible alternative pricing plans to investigate their impact on electricity cost and the number of switchings. These three pricing plans are

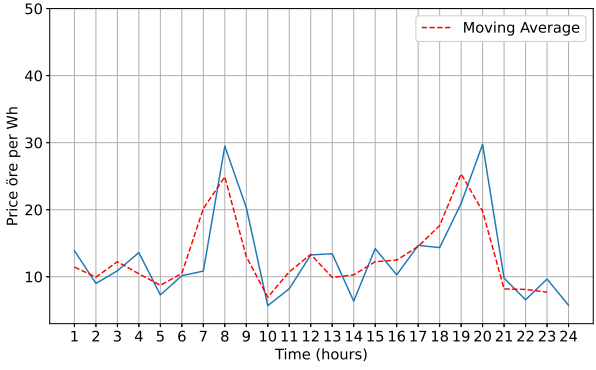
(a) Base case with standard peak hours

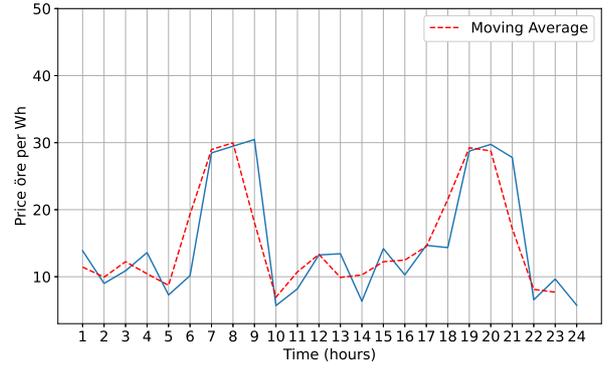
(b) Case with extended peak hours

**Figure 18:** Comparison of hourly electricity prices across two different test case.

illustrated in Fig. 20. The base policy is the same as the one we previously introduced, and the other two pricing plans stabilize price fluctuations, making it more convenient for the operator to plan their operations but with some additional cost. By solving the optimization problem for these three pricing plans for an instance with $B = 50$, $N = 25$, and $\gamma = 6$, we obtain three different Gantt charts, as shown in Fig. 21. A small instance is used here to save space for illustration. For a realistic instance with 100 batteries, the resulting Gantt chart (for the base price) is provided in the appendix (see Fig. 22). It is evident that under the base price plan (see Fig. 21c), although the electricity cost is low, the corresponding number of switchings is high. Plan 1 performs slightly better than Plan 2 in terms of switching, but it is worse in terms of electricity cost. Operators should purchase a particular pricing plan only if they believe it will result in an optimal solution with significant savings in the number of switchings, which might be desirable to reduce the cost of having staff available in the facility for a longer duration. However, it depends upon the utility metrics used by the operator to value electricity cost and manual labor.

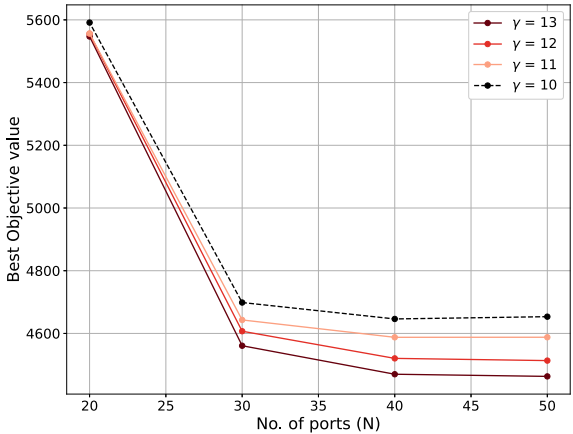
(a) Senstivity analysis for $B = 100^a$

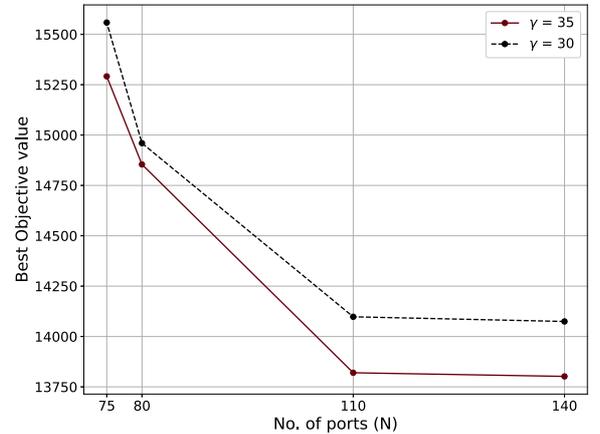
(b) Senstivity analysis for $B = 300^a$.

$^a\gamma = 9$ for $N = 50$ and $N = 10$ for $\gamma = 13$ did not result in any feasible solution

$^a\gamma = 25$ for $N = 140$ and $N = 50$ for $\gamma = 35$ did not result in any feasible solution

**Figure 19:** Senstivity analysis of variation of $\gamma$ and $N$ on the obtained objective value.

## 8. Conclusion and limitations

In this work, we present a binary linear optimization model for scheduling at a charging facility, taking micromobility depot as an example. Due to the high dimension of the problem and ineffectiveness of existing methods, we propose

Time-block enhanced Lagrangian decomposition for scheduling

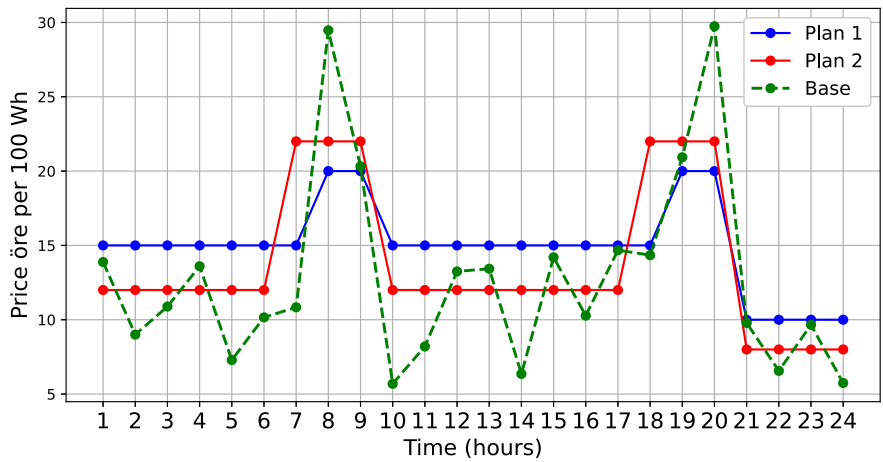

**Figure 20:** Different price plans and the base price that depend on spot prices. 100 öre = 1 SEK

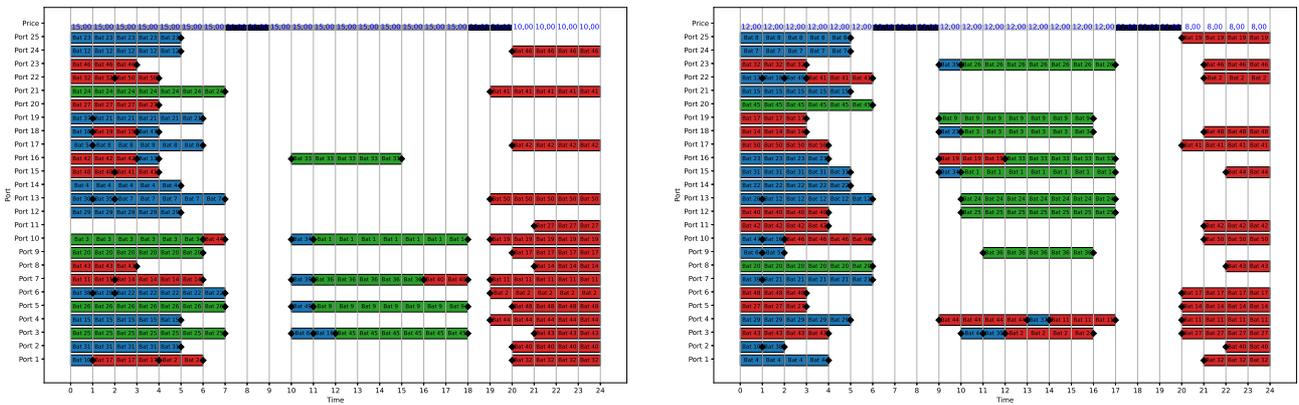

(a) Gantt chart for Plan 1. El. cost = 3165, switchings = 70.

(b) Gantt chart for Plan 2. El. cost = 2544, switchings = 77.

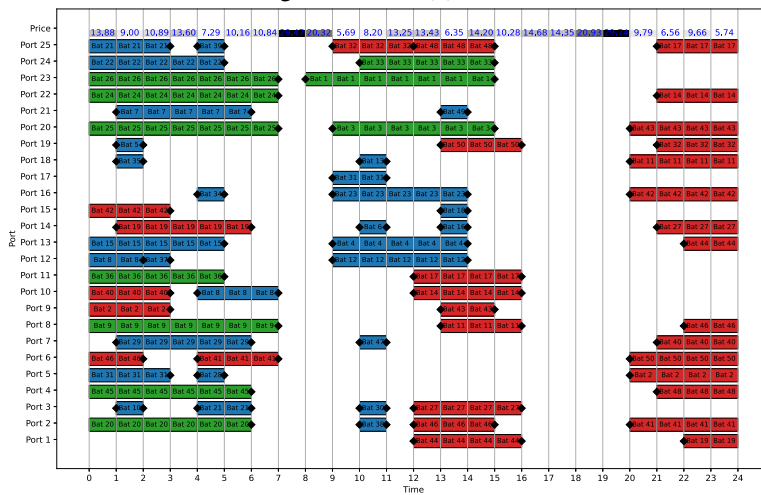

(c) Gantt chart for base plan. El. cost = 2232, switchings = 105.

**Figure 21:** Comparison of Gantt charts for $B = 50$, $N = 25$ $\gamma = 6$ using different price plans as indicated in Fig. 20

a Lagrangian decomposition approach to decompose the problem to two subproblems. We analyzed the polyhedral properties of both the subproblems. For the subproblem that posed a computational bottleneck, we suggested a new time-block formulation, which has proven to be much tighter and stronger than the original formulation. This resulted in accelerated convergence of subgradient optimization. Furthermore, we developed a new Lagrangian heuristic that utilizes the partial solutions of subproblems to fix a subset of the variables. Subsequently, a new disjunctive time-block-based bin-packing formulation is proposed, which is much more compact than simply fixing the variables in the original model. This is followed by our new ergodic-iterate-based local search method to further improve the feasible solutions obtained.

We believe that similar formulations, along with new side constraints, can describe several other scheduling applications, such as vehicle-to-grid charging and energy storage for second-life battery usage. This expansion indicates that our results can be valuable to a broader scientific community. Furthermore, to the best of our knowledge, using input parameters of a model to estimate a good initialization of Lagrangian multipliers to obtain strong starting lower bounds in subgradient optimization is not explored in-depth in current literature. This is certainly an interesting contribution that warrants further investigation.

Our proposed algorithm and polyhedral analysis can be applied to both deterministic and stochastic scheduling problems. Although our original model does not incorporate the stochasticity of, for example, demand and electricity prices, the deduced reformulations and decomposition approach can be extremely useful for solving stochastic programming problems. These problems typically involve decomposing over scenarios, where each scenario resembles our original model but with fixed parameters for each scenario.

The limitations of our method include that the time-block reformulation of the subproblem in $(\mathbf{y}, \mathbf{s})$ has a set of constraints (12b), that grows as an exponential function of the number of time periods. Therefore, if the number of time periods increases drastically, the effectiveness of our algorithm needs to be further analyzed and adjusted for other applications. Moreover, our proposal to use smaller instances to learn the relationship between Lagrangian multipliers and costs may theoretically never converge, necessitating further theoretical investigation to understand the limitations and use cases. Finally, a good extension of this work can be incorporating uncertainty in demand and electricity prices utilizing stochastic programming approach.

## CRediT authorship contribution statement

Sunney Fotedar: Writing – original draft, Methodology, Formal analysis, Visualization, Validation. Jiaming Wu: Conceptualization, Writing – review and editing, Methodology, Formal analysis, Project administration, Funding acquisition. Balázs Kulcsár: Writing – review and editing, Methodology, Funding acquisition. Rebecka Jörnsten: Writing – review & editing, Methodology, Funding acquisition.

## Acknowledgments

This work was supported by the Swedish Energy Agency under project FEAT (P2022-00404), and the Energy Area of Advance at Chalmers University of Technology under project "Energy efficiency analysis for co-existing shared and electric micromobility.".

## Appendix 1

**Proposition** ($P^\mathbf{x}$ is an integer polyderon (Prop:4.3)). *The polyhedron $P^\mathbf{x}$ has integer extreme points and corresponding subproblem in $\mathbf{x}$ is equivalent to solving a min-cost network flow problem*

*Proof.* Let $G = (V, A)$ be a directed graph where:

- $V = \{j \mid j \in \mathcal{B}\} \cup \{(j,t) \mid j \in \mathcal{B}, t \in \mathcal{T}\} \cup \{(k,t) \mid k \in \mathcal{N}, t \in \mathcal{T}\} \cup \{\text{sink}\}$

- $A = \{(j,(j,t)) \mid j \in \mathcal{B}, t \in \mathcal{T}\} \cup \{((j,t),(k,t)) \mid j \in \mathcal{B}, k \in \mathcal{N}, t \in \mathcal{T}\} \cup \{((k,t),\text{sink}) \mid k \in \mathcal{N}, t \in \mathcal{T}\}$

Each battery $j \in \mathcal{B}_\ell, \ell \in \mathcal{L}$ has a source node which has an outgoing arc to nodes corresponding to $(j,t), t \in \mathcal{T}_\ell$. Each node corresponding to $(j,t)$ has outgoing arcs to $(k,t), k \in \mathcal{N}$. Finally, each $(k,t)$ node is linked to the same sink node. Each source node has supply equal to the demand $p_j, j \in \mathcal{B}_\ell$. Each arc has a capacity limitation of one. The cost is associated only with the arc connecting $(j,t)$ nodes to $(k,t)$ nodes equal to $c_t$. Hence, if we solve this min-cost network flow problem we also solve the problem corresponding to the subproblem $P^\mathbf{x}$ □

**Proposition** (Equivalence of $F_r$ and $F_{tb}$ (Prop: 4.5)). *The time block formulation $F_{tb}$ admits an optimal solution $(\mathbf{y}^*, \lambda^*)$ within the set $\tilde{P}^{\mathbf{y}\text{-}\lambda}$. This solution can be correspondingly mapped to $(\mathbf{y}^*, \mathbf{s}^*) \in \tilde{P}^{\mathbf{y}\text{-s}}$ and retains optimality under the formulation $F_r$.*

*Proof.* We prove that if $(\mathbf{y}^*, \lambda^*) \in \tilde{P}^{\mathbf{y}\text{-}\lambda}$ is optimal for $F_{tb}$ then a corresponding solution $(\mathbf{y}^*, \mathbf{s}^*) \in \tilde{P}^{\mathbf{y}\text{-s}}$ exists that is optimal in the formulation $F_r$. Furthermore, we also prove that the two optimal solutions have the same objective value. First, we define a model to map the solution $(\mathbf{y}^*, \lambda^*) \in \tilde{P}^{\mathbf{y}\text{-}\lambda}$ to $(\mathbf{y}^*, \mathbf{s}^*)$ that is feasible in $\tilde{P}^{\mathbf{y}\text{-s}}$. Note that the switching variables are the same.

$$\min \sum_{\ell \in \mathcal{L}} \sum_{t \in \mathcal{T}_\ell} \sum_{j \in \mathcal{B}_\ell} \sum_{k \in \mathcal{N}} \bar{c}_{jkt} s_{jkt} + w \sum_{k \in \mathcal{N}} \sum_{t \in \mathcal{T} \setminus \{T\}} y_{kt}, \tag{19a}$$

$$\text{s.t.} \quad \pi_{jkt_1 t_2} \geq \frac{1}{(t_2 - t_1)} \sum_{t=t_1+1}^{t_2} s_{jkt}, \quad \forall j \in \mathcal{B}_\ell, \ell, k, (t_1, t_2) \in \mathcal{T}_<^k, \tag{19b}$$

$$\sum_{\ell \in \mathcal{L}} \sum_{j \in \mathcal{B}_\ell} \pi_{jkt_1 t_2} = \lambda_{kt_1 t_2}^*, \quad k \in \mathcal{N}, (t_1, t_2) \in \mathcal{T}_<^k, \tag{19c}$$

$$y_{kt} = y_{kt}^*, \quad k \in \mathcal{N}, t \in \mathcal{T} \setminus \{T\}, \tag{19d}$$

$$s_{jkt} \in \{0, 1\}, \quad j \in \mathcal{B}_\ell, k \in \mathcal{N}, t \in \mathcal{T} \setminus \{T\}, \ell, \tag{19e}$$

$$\pi_{jkt_1 t_2} \in \{0, 1\}, \quad j \in \mathcal{B}_\ell, k \in \mathcal{N}, (t_1, t_2) \in \mathcal{T}_<^k, \ell. \tag{19f}$$

The artificial variable $\pi_{jkt_1t_2}$ is one if the battery $j$ is used in a time block $(t_1, t_2)$ which for an optimal solution in $\tilde{P}^{\mathbf{y}\text{-}\lambda}$ we know must be with least cost $j \in \arg\min_{j \in \mathcal{B}_\ell, \ell \in \mathcal{L}} \sum_{t=t_1+1}^{t_2} \bar{c}_{jkt}$, where $\bar{c}_{jkt} < 0$. The constraint (19c) ensures that if $\lambda_{kt_1t_2}^* = 1$ then for exactly one of the $j \in \mathcal{B}_\ell$, we have $\pi_{jkt_1t_2} = 1$. This subsequently imposes a constraint $\frac{1}{(t_2-t_1)} \sum_{t=t_1+1}^{t_2} s_{jkt} \leq 1$ for the given $k$ and $(t_1, t_2)$. Since $s_{jkt}$ is binary it implies all $s_{jkt} = 1, t \in \{t_1 + 1, \ldots, t_2\}$ if $j$ has the least cost coefficient as described. We can sum over $j$ and $\ell$ (19b) and replace the l.h.s. with $\lambda_{kt_1t_2}^*$ (see (19c)) we get the following:

$$\lambda_{kt_1 t_2}^* \geq \frac{1}{(t_2 - t_1)} \sum_{\ell \in \mathcal{L}} \sum_{j \in \mathcal{B}_\ell} \sum_{t=t_1+1}^{t_2} s_{jkt}, \quad k \in \mathcal{N}, (t_1, t_2) \in \mathcal{T}_<^k. \tag{20}$$

Now we need to show that an optimal solution $(\bar{\mathbf{y}}, \bar{\mathbf{s}})$ to model (19) is also feasible in $\tilde{P}^{\mathbf{y}\text{-s}}$ and has the same objective value and if $(\mathbf{y}^*, \lambda^*) \in \tilde{P}^{\mathbf{y}\text{-}\lambda}$ is optimal in $F_{tb}$ then corresponding solution to (19) should be optimal for the formulation $F_r$. First, we show that optimal solution to (19) $(\bar{\mathbf{y}}, \bar{\mathbf{s}})$ is indeed feasible in $\tilde{P}^{\mathbf{y}\text{-s}}$. Since (20) is satisfied by $\bar{\mathbf{s}}$ then following holds:

$$\lambda_{kt_1 t_2}^* (t_2 - t_1) \geq \sum_{t=t_1+1}^{t_2} \left\{ \sum_{\ell \in \mathcal{L}} \sum_{j \in \mathcal{B}_\ell} \bar{s}_{jkt} \right\} \implies 1 \geq \sum_{\ell \in \mathcal{L}} \sum_{j \in \mathcal{B}_\ell} \bar{s}_{jkt}$$

The last inequality (which is one of the constraints in $\tilde{P}^{\mathbf{y}\text{-s}}$) comes from the fact that for $\lambda_{kt_1t_2}^* = 1$, we have $\bar{s}_{\bar{j}kt} = 1, t \in \{t_1 + 1, \ldots, t_2\}$ for one of the $\bar{j} \in \mathcal{B}_\ell, \ell \in \mathcal{L}$ as $\pi_{\bar{j}kt_1t_2} = 1$ (see (19c)). Next we ensure the switching variables $\bar{\mathbf{y}}$ also satisfy respective constraints in $\tilde{P}^{\mathbf{y}\text{-s}}$. Note that $\mathbf{y}^* = \bar{\mathbf{y}}$ as per model(19). Hence, following holds:

$$y_{kt}^* = \bar{y}_{kt} \geq \sum_{t'=t+1}^{T} \lambda_{ktt'}^* \implies \bar{y}_{kt} \geq \lambda_{ktt+1}^* \overset{(20)}{\implies} \bar{y}_{kt} \geq \sum_{\ell \in \mathcal{L}} \sum_{j \in \mathcal{B}_\ell} \bar{s}_{jkt+1} \implies \bar{y}_{kt} \geq \sum_{\ell \in \mathcal{L}} \sum_{j \in \mathcal{B}_\ell} \bar{s}_{jkt+1} - \sum_{\ell \in \mathcal{L}} \sum_{j \in \mathcal{B}_\ell} \bar{s}_{jkt} \tag{21}$$

The constraint (2g) is also deduced in a similar way. Lastly, constraint (12e) (the constraint corresponding $\gamma$) is exactly the same in $\tilde{P}^{\mathbf{y}\text{-}\mathbf{s}}$ and $\tilde{P}^{\mathbf{y}\text{-}\lambda}$. Hence, we have an optimal solution $\bar{\mathbf{y}}, \bar{\mathbf{s}}$ to model (19) that is feasible in $\tilde{P}^{\mathbf{y}\text{-}\mathbf{s}}$. It is easy to see that the objective function's value is the same due to how $d^k_{t_1 t_2}$ is defined in $F_{\text{tb}}$. We utilize the Prop. 4.4 to show that the obtained feasible solution $(\bar{\mathbf{y}}, \bar{\mathbf{s}}) \in \tilde{P}^{\mathbf{y}\text{-}\mathbf{s}}$ is also optimal for the formulation $F_{\text{r}}$. The proposition follows as given an optimal solution for time-block formulation $(\mathbf{y}^*, \lambda^*)$ we obtained the optimal solution to the formulation $F_{\text{r}}$. □

**Proposition** (Tightness of the new formulation (Proposition 4.6) ). *Show that* $\text{proj}_{\mathbf{y},\mathbf{s}}(\tilde{P}^{\mathbf{y}\text{-}\lambda}_{\text{LP}}) \subset \tilde{P}^{\mathbf{y}\text{-}\mathbf{s}}_{\text{LP}}$.

*Proof.* Let $(\bar{\mathbf{y}}, \bar{\lambda}) \in \tilde{P}^{\mathbf{y}\text{-}\lambda}_{\text{LP}}$, we need to first define a mapping $M$ of $\bar{\lambda}$ to $\bar{\mathbf{s}}$. Then we show that the resulting $(\bar{\mathbf{y}}, \bar{\mathbf{s}}) \in \tilde{P}^{\mathbf{y}\text{-}\mathbf{s}}_{\text{LP}}$. This ensure that $\text{proj}_{\mathbf{y},\mathbf{s}}(\tilde{P}^{\mathbf{y}\text{-}\lambda}_{\text{LP}}) \subseteq \tilde{P}^{\mathbf{y}\text{-}\mathbf{s}}_{\text{LP}}$ and for strict subset we provide an example when a feasible solution $(\hat{\mathbf{y}}, \hat{\mathbf{s}}) \in \tilde{P}^{\mathbf{y}\text{-}\mathbf{s}}_{\text{LP}}$ but $(\hat{\mathbf{y}}, \hat{\mathbf{s}}) \notin \text{proj}_{(\mathbf{y},\mathbf{s})}(\tilde{P}^{\mathbf{y}\text{-}\lambda}_{\text{LP}})$. Given $(\bar{\mathbf{y}}, \bar{\lambda}) \in \tilde{P}^{\mathbf{y}\text{-}\lambda}_{\text{LP}}$, we obtain corresponding mapping as discussed in previous proof in (19). Hence, for a given $(\bar{\mathbf{y}}, \bar{\lambda}) \in \tilde{P}^{\mathbf{y}\text{-}\lambda}$ we can obtain a feasible $(\bar{\mathbf{y}}, \bar{\mathbf{s}}) \in \tilde{P}^{\mathbf{y}\text{-}\mathbf{s}}$ as follows[9]

$$\bar{\mathbf{s}} \in \arg\min \left\{ \sum_{\ell \in \mathcal{L}} \sum_{t \in \mathcal{T}_\ell} \sum_{j \in \mathcal{B}_\ell} \sum_{k \in \mathcal{N}} \bar{c}_{jkt} s_{jkt} \,\middle|\, \bar{\lambda}_{k t_1 t_2} \geq \frac{1}{(t_2 - t_1)} \sum_{\ell \in \mathcal{L}} \sum_{j \in \mathcal{B}_\ell} \sum_{t=t_1+1}^{t_2} s_{jkt}, \quad k \in \mathcal{N}, (t_1, t_2) \in \mathcal{T}^k_< \right\}. \quad (22)$$

The existence of such $\bar{\mathbf{s}}$ can be deduced by reformulating as follows:

$$\bar{\lambda}_{k t_1 t_2}(t_2 - t_1) \geq \sum_{\ell \in \mathcal{L}} \sum_{j \in \mathcal{B}_\ell} \sum_{t=t_1+1}^{t_2} s_{jkt}, \quad k \in \mathcal{N}, (t_1, t_2) \in \mathcal{T}^k_<. \quad (23)$$

Note that if there does not exist $\bar{\mathbf{s}}$ that satisfies this constraint then it can be shown that $\tilde{P}^{\mathbf{y}\text{-}\mathbf{s}} \in \emptyset$ as follows. The constraint (8d) of $\tilde{P}^{\mathbf{y}\text{-}\mathbf{s}}$ can be aggregated to get $\sum_{\ell \in \mathcal{L}} \sum_{j \in \mathcal{B}_\ell} \sum_{t=t_1+1}^{t_2} s_{jkt} \leq (t_2 - t_1)$ and since $\lambda_{k, t_1, t_2} \in [0, 1]$ it is not possible that no $\mathbf{s}$ satisfies (23) unless $\tilde{P}^{\mathbf{y}\text{-}\mathbf{s}} \in \emptyset$. Now we show that $\bar{\mathbf{s}}$ satisfies constraints in $\tilde{P}^{\mathbf{y}\text{-}\mathbf{s}}$. By the definition of $\bar{\mathbf{s}}$ and since $\bar{\lambda}_{k t_1 t_2} \in [0, 1]$ it satisfies constraint (8d). Now we need to show that $\bar{\mathbf{s}}$ satisfies the switching constraints (8b)-(8c) for a given $\bar{\mathbf{y}}$. First we know that $\bar{y}_{k t_1} \geq \sum_{t=t_1+1}^T \lambda_{k, t_1, t}$. Now by using (23) we can deduce following:

$$\bar{y}_{k t_1} \geq \sum_{t=t_1+1}^T \left\{ \frac{1}{(t - t_1)} \sum_{t'=t_1+1}^t \sum_{\ell \in \mathcal{L}} \sum_{j \in \mathcal{B}_\ell} \bar{s}_{jkt'} \right\} \implies \bar{y}_{k t_1} \geq \sum_{\ell \in \mathcal{L}} \sum_{j \in \mathcal{B}_\ell} \bar{s}_{j k t_1 + 1}$$

The last inequality also implies (8b). Lastly, we deduce that $\bar{\mathbf{s}}$ satisfy (8c). We start with constraint $\bar{y}_{k t_2} \geq \sum_{t=0}^{t_2-1} \lambda_{k t t_2}, k \in \mathcal{N}, t_2 \in \mathcal{T} \setminus \{T\}$. Now using (23) we deduce following:

$$\bar{y}_{k t_2} \geq \sum_{t=0}^{t_2-1} \left\{ \frac{1}{(t_2 - t)} \sum_{t'=t+1}^{t_2} \sum_{\ell \in \mathcal{L}} \sum_{j \in \mathcal{B}_\ell} s_{jkt'} \right\} \implies \bar{y}_{k t_2} \geq \sum_{\ell \in \mathcal{L}} \sum_{j \in \mathcal{B}_\ell} \bar{s}_{j k t_2} \implies \bar{y}_{k t_2} \geq \bar{s}_{j k t_2}$$

The last inequality implies that $\bar{\mathbf{s}}$ also satisfies (8c). Furthermore $\bar{\mathbf{y}}$ already satisfies the constraint on total number of switching in a time period. Hence, any solution $(\bar{\mathbf{y}}, \bar{\lambda}) \in \tilde{P}^{\mathbf{y}\text{-}\lambda}$ can be projected onto $(\mathbf{y}, \mathbf{s})$ such that $\text{proj}_{\mathbf{y},\mathbf{s}}(\tilde{P}^{\mathbf{y}\text{-}\lambda}) \subseteq \tilde{P}^{\mathbf{y}\text{-}\mathbf{s}}$. To prove strict inequality we provide an example where a solution $(\hat{\mathbf{y}}, \hat{\mathbf{s}}) \in \tilde{P}^{\mathbf{y}\text{-}\mathbf{s}}$ but $(\hat{\mathbf{y}}, \hat{\mathbf{s}}) \notin \text{proj}_{\mathbf{y},\mathbf{s}}(\tilde{P}^{\mathbf{y}\text{-}\lambda})$. Let $N = 1, T = 5, \mathcal{B}_1 = \mathcal{B} = \{1, 2, 3, 4\}, L = 1$ (similar to the example in Fig. 2). Next let us consider a feasible solution $(\hat{\mathbf{y}}, \hat{\mathbf{s}}) \in \tilde{P}^{\mathbf{y}\text{-}\mathbf{s}}_{\text{LP}}$. We have $\hat{s}_{j1t} = 0, j \in \{1, 2, 3\}, t \in \{1, 2, 3, 4, 5\}, \hat{s}_{411} = \hat{s}_{412} = \hat{s}_{413} = 0.333, \hat{s}_{414} = 0.1$, $\hat{s}_{415} = 0.9$ and $y_{1t} = 0, t = \{1, 2, 4\}$. Due to the constraint (8c) we have $\hat{y}_{13} \geq \hat{s}_{413} - \hat{s}_{414} = 0.333 - 0.1 = 0.233$. Hence, we can have a feasible solution $\hat{y}_{13} = 0.233$. Now from (23), we get $\hat{\lambda}_{k35}(2) \geq 1 \implies \hat{\lambda}_{k35} \geq 0.5$. Now using the constraint (12c) we get $\hat{y}_{13} \geq \hat{\lambda}_{34} + \hat{\lambda}_{35} \implies \hat{y}_{13} \geq 0.5$ which is clearly not true in the considered case. The proposition follows. □

---

[9] Note that the mapping is slightly different than in the previous proof as we have the relaxed set in this proof.

# Appendix 2

**Optimization model for Local search** For a given $\sigma$, $\mathcal{N}_\sigma(\tilde{x}^i)$ and $\mathcal{B}_\ell(\mathbf{x}^{\text{best}}, \mathcal{N}_\sigma(\tilde{x}^i))$ we consider local search optimization problem as follows[10].:

$$z^*(\sigma, \mathbf{x}^{\text{best}}, \tilde{\mathbf{x}}^i) := \min\left(\sum_{\ell \in \mathcal{L}} \sum_{t \in \mathcal{T}_\ell} c_t \left(\sum_{j \in \mathcal{B}_{\ell,\sigma}} \sum_{k \in \mathcal{N}_\sigma} x_{jkt}\right) + w \sum_{k \in \mathcal{N}_\sigma} \sum_{t \in \mathcal{T} \setminus T} y_{kt}\right), \tag{24a}$$

$$\text{s.t.} \sum_{t \in \mathcal{T}_\ell} \left(\sum_{k \in \mathcal{N}_\sigma} x_{jkt}\right) = \sum_{t \in \mathcal{T}} \left(\sum_{k \in \mathcal{N}_\sigma} x_{jkt}^{\text{best}}\right), \qquad j \in \mathcal{B}_{\ell,\sigma}, \ell \in \mathcal{L}, \tag{24b}$$

$$\sum_{\ell \in \mathcal{L}} \sum_{j \in \mathcal{B}_{\ell,\sigma}} x_{jkt} \leq 1, \qquad k \in \mathcal{N}_\sigma, t \in \mathcal{T}, \tag{24c}$$

$$\sum_{k \in \mathcal{N}_\sigma} x_{jkt} + \sum_{k' \in \mathcal{N} \setminus \mathcal{N}_\sigma} x_{jk't}^{\text{best}} \leq 1, \qquad j \in \mathcal{B}_{\ell,\sigma} t \in \mathcal{T}_\ell, \ell \in \mathcal{L}, \tag{24d}$$

$$y_{kt} \geq \left(\sum_{\ell \in \mathcal{L}} \sum_{j \in \mathcal{B}_{\ell,\sigma}} x_{jkt+1} - \sum_{\ell \in \mathcal{L}} \sum_{j \in \mathcal{B}_{\ell,\sigma}} x_{jkt}\right), \qquad k \in \mathcal{N}_\sigma, t \in \mathcal{T} \setminus \{T\}, \tag{24e}$$

$$y_{kt} \geq (x_{jkt} - x_{jkt+1}), \qquad j \in \mathcal{B}_{\ell,\sigma}, k \in \mathcal{N}_\sigma, t \in \mathcal{T}_\ell \setminus \{T\}, \ell \in \mathcal{L}, \tag{24f}$$

$$\sum_{k \in \mathcal{N}_\sigma} y_{kt} + \sum_{k \in \mathcal{N} \setminus \mathcal{N}_\sigma} y_{kt}^{\text{best}} \leq \gamma, \qquad t \in \mathcal{T} \setminus \{T\}, \tag{24g}$$

$$x_{jkt} \in \{0,1\}, \qquad j \in \mathcal{B}_{\ell,\sigma}, k \in \mathcal{N}_\sigma, t \in \mathcal{T}_\ell, \ell \in \mathcal{L}, \tag{24h}$$

$$y_{kt} \in \{0,1\}, \qquad k \in \mathcal{N}_\sigma, t \in \mathcal{T} \setminus \{T\}. \tag{24i}$$

**Ergodic iterate updates** To avoid storing Lagrangian subproblem solution $\mathbf{x}^*(\boldsymbol{\mu}^i)$ for all the iteration (Gustavsson et al., 2014) suggests:

$$\tilde{\mathbf{x}}^1 = \mathbf{x}^*(\boldsymbol{\mu}^i), \qquad \tilde{\mathbf{x}}^i := \frac{\tilde{\alpha}_{i-1}}{\tilde{\alpha}_i}\tilde{\mathbf{x}}^{i-1} + \frac{\alpha_{i-1}}{\tilde{\alpha}_i}\mathbf{x}^*(\boldsymbol{\mu}^{i-1}), \quad i = 2,3,\ldots. \quad \tilde{\alpha}_1 := \alpha_0, \quad \tilde{\alpha}_i := \tilde{\alpha}_{i-1} + \alpha_{i-1}, i = 1,2,\ldots.$$
$$\tag{25}$$

---

[10]For simplification of notation we use $\mathcal{N}_\sigma$ instead of $\mathcal{N}_\sigma(\tilde{x}^i)$ and $\mathcal{B}_{\ell,\sigma}$ instead of $\mathcal{B}_\ell(\mathbf{x}^{\text{best}}, \mathcal{N}_\sigma(\tilde{x}^i))$

Time-block enhanced Lagrangian decomposition for scheduling

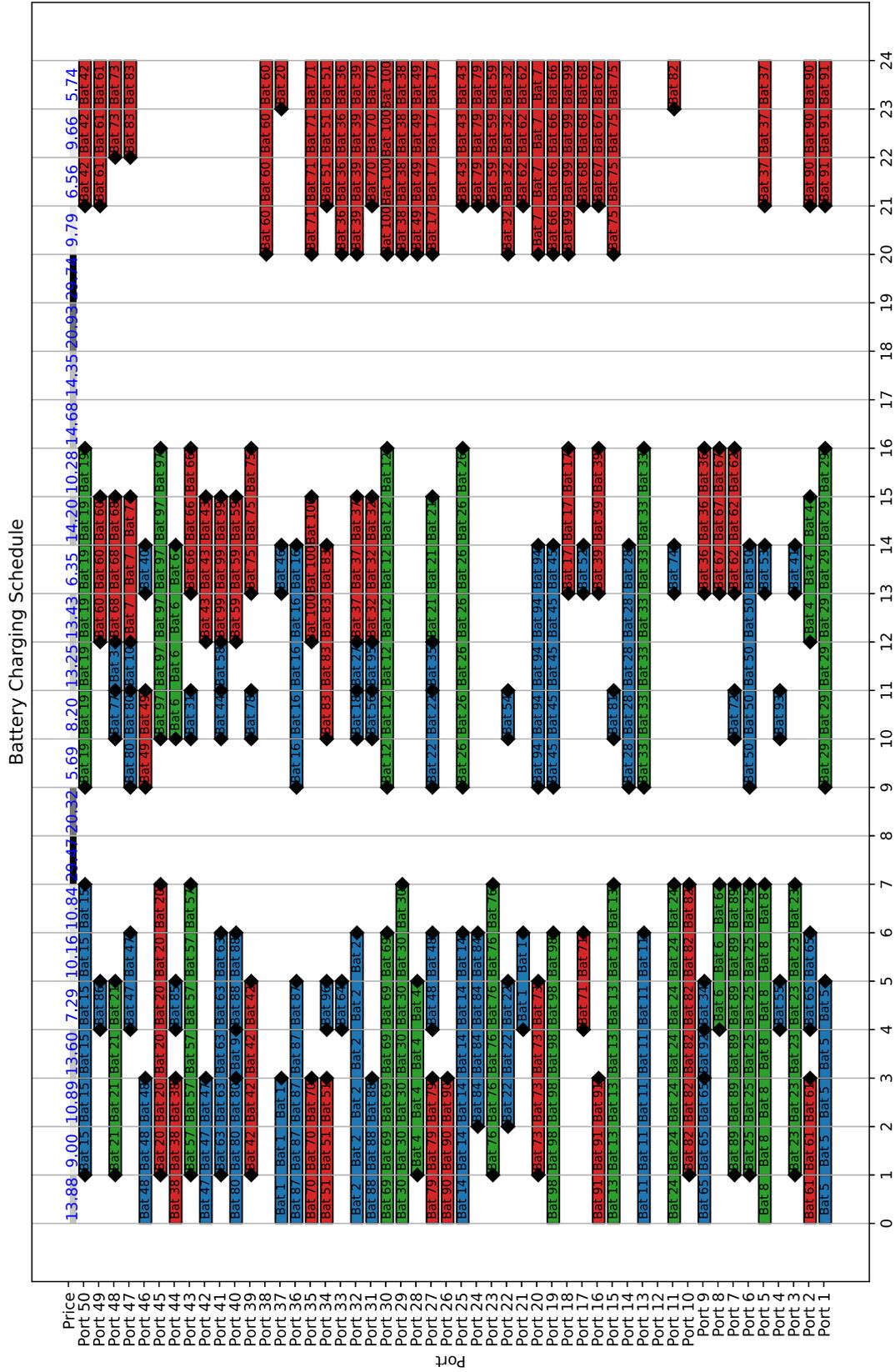

**Figure 22:** Gantt chart for the instance with $B = 100, N = 50, \gamma = 13$, and electricity price is mentioned at the top. The batteries are color coded blue: $\mathcal{B}_1$, green: $\mathcal{B}_2$, and red: $\mathcal{B}_3$.

Fotedar et al.: *Preprint submitted to Elsevier* Page 37 of 39

## Data availability

Data and code can be found in: https://github.com/SunneyF/Lagrangian-Decomposition-for-Battery-Charging.

NACTO, 2024. A micromobility record: 157 million trips on bike share and scooter share in 2023. URL: https://nacto.org/2024/07/22/a-micromobility-record-157-million-trips-on-bike-and-scooter-share-in-2023/.

Nayak, D.S., Misra, S., 2024. An operational scheduling framework for electric vehicle battery swapping station under demand uncertainty. Energy 290, 130219. doi:mgqb.

Nurre, S.G., Bent, R., Pan, F., Sharkey, T.C., 2014. Managing operations of plug-in hybrid electric vehicle (phev) exchange stations for use with a smart grid. Energy Policy 67, 364–377. doi:mgvh.

Pelletier, S., Jabali, O., Laporte, G., 2018. Charge scheduling for electric freight vehicles. Transportation Research Part B: Methodological 115, 246–269. doi:gd7vms.

Polyak, B., 1969. Minimization of unsmooth functionals. USSR Computational Mathematics and Mathematical Physics 9, 14–29. doi:cthfpb.

Raeesi, R., Zografos, K.G., 2022. Coordinated routing of electric commercial vehicles with intra-route recharging and en-route battery swapping. European Journal of Operational Research 301, 82–109. doi:mhcw.

Raviv, T., 2012. The battery switching station scheduling problem. Operations Research Letters 40, 546–550. doi:f4hv7b.

Sagastizábal, C., 2012. Divide to conquer: decomposition methods for energy optimization. Mathematical Programming 134, 187–222. doi:w6z.

Severson, K.A., Attia, P.M., Jin, N., Perkins, N., Jiang, B., Yang, Z., Chen, M.H., Aykol, M., Herring, P.K., Fraggedakis, D., et al., 2019. Data-driven prediction of battery cycle life before capacity degradation. Nature Energy 4, 383–391.

Sherali, H.D., Ulular, O., 1989. A primal-dual conjugate subgradient algorithm for specially structured linear and convex programming problems. Applied Mathematics & Optimization 20, 193–221. doi:cmgf25.

Sun, B., Tan, X., Tsang, D.H.K., 2018. Optimal charging operation of battery swapping and charging stations with qos guarantee. IEEE Transactions on Smart Grid 9, 4689–4701. doi:gd6ttz.

Tan, X., Qu, G., Sun, B., Li, N., Tsang, D.H.K., 2019. Optimal scheduling of battery charging station serving electric vehicles based on battery swapping. IEEE Transactions on Smart Grid 10, 1372–1384. doi:gf8vn8.

VOI Technology, 2024. Safety standards of battery charging at voi technology's warehouse. LinkedIn. URL: https://shorturl.at/ckBK1.

Wolsey, L.A., 2001. Integer Programming. volume 1. Wiley.

Zhan, W., Wang, Z., Zhang, L., Liu, P., Cui, D., Dorrell, D.G., 2022. A review of siting, sizing, optimal scheduling, and cost-benefit analysis for battery swapping stations. Energy 258, 124723. doi:pkq8.

Zhao, J., Wu, J., Fotedar, S., Li, Z., Liu, P., 2024. Fleet availability analysis and prediction for shared e-scooters: An energy perspective. Transportation Research Part D: Transport and Environment 136, 104425. doi:nj44.

Zhao, Y., Larsson, T., Rönnberg, E., Pardalos, P.M., 2018. The fixed charge transportation problem: a strong formulation based on Lagrangian decomposition and column generation. Journal of Global Optimization 72, 517–538. doi:gd8hf7.

Özarık, S.S., Lurkin, V., Veelenturf, L.P., Van Woensel, T., Laporte, G., 2023. An adaptive large neighborhood search heuristic for last-mile deliveries under stochastic customer availability and multiple visits. Transportation Research Part B: Methodological 170, 194–220. doi:grzbk6.